\documentclass[12pt]{iopart}

\usepackage{graphicx}
\usepackage{float}
\usepackage{cases}
\usepackage{multirow}
\usepackage{color}

\usepackage{iopams}  

\newtheorem{theorem}{Theorem}[section]
\newtheorem{lemma}{Lemma}[section]
\newtheorem{corollary}{Corollary}[section]
\newtheorem{remark}{Remark}[section]
\newtheorem{definition}{Definition}[section]

\newcommand{\bproof}{{\noindent \bf Proof }}
             \newcommand{\eproof}{\hfill$\square$}

\begin{document}

\title[Shape Reconstruction in Linear Elasticity]{Shape Reconstruction in Linear Elasticity: 
\\
Standard and Linearized Monotonicity Method}

\author{Sarah Eberle  and Bastian Harrach}

\address{ Institute of Mathematics,
Goethe-University Frankfurt, Frankfurt am Main, Germany}
\ead{eberle@math.uni-frankfurt.de$,$ harrach@math.uni-frankfurt.de}

\begin{abstract}
In this paper, we deal with the inverse problem of the shape reconstruction of inclusions in elastic bodies.
The main idea of this reconstruction is based on the monotonicity property of the Neumann-to-Dirichlet operator presented in a former article of the authors.
Thus, we introduce the so-called standard as well as linearized monotonicity tests in order to detect and reconstruct inclusions.
In addition, we compare these methods with each other and present several numerical
test examples.
\end{abstract}

%

\section{Introduction and problem statement}
The reconstruction of inclusions in materials by nondestructive testing becomes more and more important and opens a wide mathematical
field in inverse problems. The applications cover engineering, geoscientific and medical problems.
This is the reason, why several authors dealt with the inverse problem of linear elasticity in order to recover the Lam\'e parameters. 
For the two dimensional case, we refer the reader to  \cite{nakamura1993identification,imanuvilov2011reconstruction,
lin2017boundary}, {\color{black} and \cite{ikehata1990inversion} deals with both the two and three-dimensional case.}
Further on, in three dimensions, \cite{nakamura1995inverse, nakamura2003global} and \cite{eskin2002inverse}
gave the proof for uniqueness results for both Lam\'e coefficients under the assumption that $\mu$ is close to a positive constant. 
\cite{beretta2014lipschitz, beretta2014uniqueness} proved the uniqueness for partial data, where the Lam\'e parameters are piecewise
constant and some boundary determination results were shown in \cite{nakamura1999layer,nakamura1995inverse,lin2017boundary}.
{\color{black} In addition, for uniqueness results concerning the anisotropic case, the reader is referred to \cite{Ikehata1999} as well as to
 \cite{Ikehata1998} for an approach based on monotonicity properties, and \cite{Ikehata2006} for an approach that does not depend on
monotonicity properties.
Finally, we want to mention \cite{Ikehata2002} for the reconstruction of inclusion from boundary measurements.}
 \\
 \\
In this paper, the key issue of the shape reconstruction of inclusions is the \mbox{monotonicity} property of the corresponding Neumann-to-Dirichlet operator 
(see \cite{Tamburrino06, tamburrino2002new}).
These monotonicity properties were also applied for electrical impedance \mbox{tomography}
(see, e.g., \cite{harrach2013monotonicity}) and for linear elasticity in \cite{Eberle_Monotonicity}, which is
the basis for our current work. Our approach relies on the monotonicity
 of the Neumann-to-Dirichlet operator with respect to the  Lam\'e parameters and the techniques of localized potentials 
\cite{harrach2018helmholtz,harrach2018localizing,harrach2012simultaneous,gebauer2008localized,harrach2009uniqueness}.
\\
Thus, we start with the introduction of the problem and summarize the main results from  
\cite{Eberle_Monotonicity}. After that, we go over to the shape reconstruction itself, where we consider the standard and linearized monotonicity method.
In doing so, we present and compare numerical experiments from the aforementioned two methods.
\\
\\
We start with the introduction of the problem of interest in the following way.
We consider a bounded and connected open set $\Omega\subset \mathbb{R}^d$ ($d=2$ or $3$) with smooth boundary, occupied by an isotropic material with linear stress-strain relation, 
where $\Gamma_{\mathrm{D}}$ and $\Gamma_{\mathrm{N}}$ with $\Gamma_{\mathrm{D}}\cup\Gamma_{\mathrm{N}}=\partial\Omega$ are the corresponding Dirichlet and Neumann boundaries. Then the displacement vector 
$u:\overline{\Omega}\to\mathbb{R}^d$ satisfies the boundary value problem
{\color{black}
\noindent
\begin{equation}\left\{ \begin{array}{rcll}\label{direct_1}
\nabla\cdot\left(\lambda (\nabla\cdot u)I + 2\mu \hat{\nabla} u \right)&=&0 &\quad \mathrm{in}\,\,\Omega,\\
\left(\lambda (\nabla\cdot u)I + 2\mu \hat{\nabla} u \right) n&=&g&\quad \mathrm{on}\,\, \Gamma_{\mathrm{N}},\\
u&=&0 &\quad \mathrm{on}\,\, \Gamma_{\mathrm{D}},
\end{array}\right.
\end{equation}
\noindent
where $\mu$, $\lambda\in L^\infty_+(\Omega)$ are the Lam\'{e} parameters,
$\hat{\nabla} u=\frac{1}{2}\left(\nabla u + (\nabla u)^T\right)$ is the symmetric gradient, 
$n$ is \mbox{the normal} vector pointing outside of $\Omega$, 
$g\in L^{2}(\Gamma_{\mathrm{N}})^d$ the boundary load and $I$ the $d\times d$-identity matrix.
\\
Here, we use the following definitions:
 \begin{equation*}
 \nabla u(x)=
\left(\frac{\partial}{\partial x_1}u(x),\cdots,\frac{\partial}{\partial x_d}u(x)\right),
\\ \quad
 \mathrm{div}F(x)=
 \left( \begin{array}{cccc}
\sum \limits_{i=1}^d \frac{\partial}{\partial x_i}F_{i1}(x)\\
\vdots\\
\sum \limits_{i=1}^d \frac{\partial}{\partial x_i}F_{id}(x)
\end{array}\right)
\end{equation*}
 \noindent
 for $F:\Omega\mapsto \mathbb{R}^{d\times d}$.
\\
\noindent
The weak  formulation of the problem (\ref{direct_1}) is that $u\in\mathcal{V}$ fulfills}
\begin{eqnarray}
\label{var-direct_1}
\int_{\Omega} 2 \mu\, \hat{\nabla}u : \hat{\nabla}v  + \lambda \nabla \cdot u \,\nabla\cdot  v\,dx=\int_{\Gamma_{\mathrm{N}}}g \cdot v \,ds \quad \mathrm{ for\,all\,} v\in \mathcal{V},
\end{eqnarray}
where
\begin{eqnarray*}
\mathcal{V}:=\left\{   v\in H^1(\Omega)^d:\,  v_{|_{\Gamma_{\mathrm{D}}}}=0\right\},
\end{eqnarray*}
\noindent
{\color{black} which is defined in a similar way as, e.g., in \cite{Ciarlet}.}
\noindent
\\
\\
The existence and uniqueness of a solution to the above variational formulation (\ref{var-direct_1}) follows from 
the Lax-Milgram theorem, see e.g., in \cite{Ciarlet}.
\noindent
\\
\\
Next, we introduce the Neumann-to-Dirichlet operator $\Lambda(\lambda,\mu)$ by
\begin{eqnarray*}
\Lambda(\lambda,\mu): L^2(\Gamma_{\mathrm{N}})^d\rightarrow L^2(\Gamma_{\mathrm{N}})^d: 
\quad  g\mapsto{\color{black} {u^g_{(\lambda,\mu)}}\Big|_{\Gamma_{\mathrm{N}}}},
\end{eqnarray*}
\noindent
{\color{black}where $u^g_{(\lambda,\mu)}\in\mathcal{V}$ is the unique solution of the boundary value problem (\ref{direct_1}) for the Lam\'e
parameter $\lambda,\mu$ and the boundary load $g\in L^2(\Gamma_{\mathrm{N}})^d$.}
\\
\\
{\color{black} This operator describes the displacement $u$, where $u\in\mathcal{V}$ is the unique solution of the boundary value
problem (\ref{direct_1}) on the Neumann boundary $\Gamma_{\mathrm{N}}$ resulting
from the application of the boundary load $g\in L^{2}(\Gamma_{\mathrm{N}})^d$ on the Neumann boundary $\Gamma_{\mathrm{N}}$.
The associated bilinear form is given by
\begin{eqnarray}\label{bilinear_form}
\langle g,\Lambda(\lambda,\mu)h\rangle=
\int_{\Omega} \left( 2 \mu\, \hat{\nabla}u^{g}_{(\lambda,\mu)} : \hat{\nabla}u^{h}_{(\lambda,\mu)}  + 
\lambda \nabla \cdot u^{g}_{(\lambda,\mu)} \,\nabla\cdot  u^{h}_{(\lambda,\mu)}\right)\,dx,\label{bilinear_Lambda}
\end{eqnarray}
\noindent
where $u_{(\lambda,\mu)}^{g}$ and $u_{(\lambda,\mu)}^{h}$ are the unique solutions of the boundary value problem (\ref{direct_1}) 
with the boundary loads $g,h\in L^{2}(\Gamma_{\mathrm{N}})^d$, respectively.}
{\color{black}
\begin{corollary}\label{Lambda_compact}
The Neumann-to-Dirichlet operator $\Lambda(\lambda,\mu)$ is a compact linear and self-adjoint operator. 
\end{corollary}
\bproof
The linearity and self-adjointness directly follow from the properties of the bilinear form (\ref{bilinear_Lambda}) and the compactness via the compactness 
of the trace operator.
\eproof}
\\
\begin{remark}
The  inverse problem we consider here is the following:
\begin{center}
 Find $(\lambda,\mu)$ knowing the Neumann-to-Dirichlet operator $\Lambda(\lambda,\mu)$.
\end{center}
\end{remark}

\section{Monotonicity methods}

In this section, we introduce two monotonicity methods in order to reconstruct inclusions in elastic bodies. The first method is the standard (or non-linearized) monotonicity
method and the second the linearized monotonicity method. For both methods, we analyze the monotonicity properties for the Neumann-to-Dirichlet operator and 
formulate the corresponding monotonicity test which we apply for the realization of the numerical experiments.

{\color{black}
\subsection{Monotonicity estimates and Fr\'echet differentiability }
}
First, we summarize and present the required results concerning the monotonicity properties. 
The details and proofs can be found in \cite{Eberle_Monotonicity}. 
We start with the monotonicity estimate and the monotonicity property itself, which is the key issue for our study and will be analyzed later on in detail.
{\color{black} Further on, we want to remark that similar results as in Lemma \ref{mono} and Lemma \ref{mon_est_2} can be found in \cite{Ikehata1998}.}
\begin{lemma}[Lemma 3.1 from \cite{Eberle_Monotonicity}]
\label{mono}
Let $(\lambda_0,\mu_0), (\lambda_1,\mu_1) \in{\color{black}L^\infty_+(\Omega)\times L^\infty_+(\Omega)}$,   $g\in L^2(\Gamma_{\mathrm{N}})^d$ be an applied boundary load, and let 
$u_0:=u^{g}_{(\lambda_0,\mu_0)}\in \mathcal{V}$, $u_1:=u^{g}_{(\lambda_1,\mu_1)}\in \mathcal{V}$. Then
\begin{eqnarray}
\label{eqmono}
&\int_\Omega 2(\mu_0-\mu_1)\hat{\nabla}u_1:\hat{\nabla}u_1+(\lambda_0-\lambda_1)\nabla\cdot u_1\nabla\cdot u_1\,dx\\ \nonumber
&\quad\geq \langle g,\Lambda(\lambda_1,\mu_1)g\rangle-\langle g,\Lambda(\lambda_0,\mu_0)g\rangle\\ 
&\quad\geq \int_\Omega 2(\mu_0-\mu_1)\hat{\nabla}u_0 :\hat{\nabla}u_0+(\lambda_0-\lambda_1)\nabla\cdot u_0 \nabla\cdot u_0\,dx.
\end{eqnarray}
\end{lemma}

\begin{lemma}\label{mon_est_2}
Let $(\lambda_0,\mu_0),(\lambda_1,\mu_1)  \in {\color{black}L^\infty_+(\Omega)\times L^\infty_+(\Omega)}$,   $g\in L^2(\Gamma_{\mathrm{N}})^d$ be an applied boundary load, and let 
\mbox{$u_0:=u^{g}_{(\lambda_0,\mu_0)}\in \mathcal{V}$}, $u_1:=u^{g}_{(\lambda_1,\mu_1)}\in \mathcal{V}$. Then
\begin{eqnarray}\label{mono_lame}
&\langle g,\Lambda(\lambda_1,\mu_1)g\rangle-\langle g,\Lambda(\lambda_0,\mu_0)g\rangle\\  \nonumber
&\quad \geq 
\int_{\Omega} 2\left(\mu_1-\frac{\mu_1^2}{\mu_0}\right)\hat{\nabla}u_1:\hat{\nabla}u_1\,dx
+\int_{\Omega}\left(\lambda_1-\frac{\lambda_1^2}{\lambda_0}\right)\nabla\cdot u_1 \nabla \cdot u_1\,dx\\
&\quad=
\int_{\Omega} 2\frac{\mu_1}{\mu_0}\left(\mu_0-\mu_1\right)\hat{\nabla}u_1:\hat{\nabla}u_1\,dx
+\int_{\Omega}\frac{\lambda_1}{\lambda_0}\left(\lambda_0-\lambda_1\right) \nabla\cdot u_1 \nabla \cdot u_1\,dx.
\end{eqnarray}
\end{lemma}

\bproof
We start with a result shown in the {\color{black}proof} of Lemma 3.1 in \cite{Eberle_Monotonicity}:
\begin{eqnarray*}
&\langle g,\Lambda(\lambda_1,\mu_1)g\rangle-\langle g,\Lambda(\lambda_0,\mu_0)g\rangle\\
&\quad=\int_\Omega 2\mu_1\hat{\nabla}(u_1-u_0):\hat{\nabla}(u_1-u_0)+\lambda_1\nabla\cdot(u_1-u_0)\nabla\cdot (u_1-u_0)\\
&\qquad+2(\mu_0-\mu_1)\hat{\nabla}u_0:\hat{\nabla}u_0+(\lambda_0-\lambda_1)\nabla\cdot u_0\nabla\cdot u_0\,dx.
\end{eqnarray*}
\noindent
Based on this, we are led to
\noindent
\begin{eqnarray*}
&\langle g,\Lambda(\lambda_1,\mu_1)g\rangle-\langle g,\Lambda(\lambda_0,\mu_0)g\rangle\\
&\quad= \int_\Omega 2\mu_1\hat{\nabla}u_1:\hat{\nabla}u_1- 4\mu_1 \hat{\nabla}u_1:\hat{\nabla}u_0+2\mu_1\hat{\nabla}u_0:\hat{\nabla}u_0+\lambda_1\nabla\cdot u_1 \nabla\cdot u_1\\
&\qquad-2\lambda_1 \nabla\cdot u_1 \nabla\cdot u_0 +\lambda_1\nabla\cdot u_0 \nabla\cdot u_0+2\mu_0\hat{\nabla}u_0:\hat{\nabla}u_0+\lambda_0\nabla\cdot u_0\nabla\cdot u_0\\
&\qquad-2\mu_1\hat{\nabla}u_0:\hat{\nabla} u_0-\lambda_1 \nabla\cdot u_0 \nabla¸\cdot u_0\,dx\\
&\quad = \int_\Omega 2\mu_1\left\Vert\hat{\nabla}u_1\right\Vert_F^2- 4\mu_1 \hat{\nabla}u_1:\hat{\nabla}u_0+\lambda_1\left\vert\nabla\cdot u_1 \right\vert^2 -2\lambda_1 \nabla\cdot u_1 \nabla\cdot u_0 \\
&\qquad+2\mu_0\left\Vert\hat{\nabla}u_0\right\Vert_F^2+\lambda_0\left\vert\nabla\cdot u_0\right\vert^2\,dx\\
&\quad=\int_\Omega 2\mu_0\left\Vert \hat{\nabla}u_0-\frac{\mu_1}{\mu_0}\hat{\nabla}u_1\right\Vert_F^2+\lambda_0\left\vert\nabla\cdot u_0-\frac{\lambda_1}{\lambda_0}\nabla\cdot u_1\right\vert^2\\
&\qquad+\left(2\mu_1-\frac{2\mu_1^2}{\mu_0}\right)\left\Vert \hat{\nabla}u_1\right\Vert_F^2
+\left(\lambda_1-\frac{\lambda_1^2}{\lambda_0}\right)\left\vert \nabla\cdot u_1\right\vert^2\,dx\\
&\quad\geq\int_{\Omega}2\left(\mu_1-\frac{\mu_1^2}{\mu_0}\right)\left\Vert \hat{\nabla}u_1\right\Vert_F^2+\left(\lambda_1-\frac{\lambda_1^2}{\lambda_0}\right)\left\vert \nabla\cdot u_1\right\vert^2\,dx.
\end{eqnarray*}
\noindent
\eproof
\\
\noindent
{\color{black}
Both, Lemma \ref{mono} and Lemma \ref{mon_est_2} lead}
directly to
\begin{corollary}[Corollary 3.2 from \cite{Eberle_Monotonicity}]\label{monotonicity}
For $(\lambda_0,\mu_0), (\lambda_1,\mu_1) \in {\color{black}L^\infty_+(\Omega)\times L^\infty_+(\Omega)}$
\begin{eqnarray*}
\lambda_0\leq \lambda_1 \,\,\textnormal{\it and }\,\, \mu_0\leq \mu_1 \,\, \textnormal{\it  implies } \,\, \Lambda(\lambda_0,\mu_0)\geq \Lambda(\lambda_1,\mu_1).
\end{eqnarray*}
\end{corollary}
\noindent
\\
{\color{black}
%
%
\noindent
Next, we are concerned with the Fr\'echet differentiability of the Neumann-to-Dirichlet operator $\Lambda$.

\begin{lemma}
The Neumann-to-Dirichlet operator $\Lambda:L^\infty_+(\Omega)\times L^\infty_+(\Omega)\to \mathcal{L}(L^2(\Gamma_{\mathrm{N}})^d)$ is 
Fr\'echet differentiable.
The Fr\'echet derivative
\begin{eqnarray*}
\Lambda'(\lambda,\mu):  L^\infty(\Omega)\times L^\infty(\Omega)\rightarrow  \mathcal{L}(L^2(\Gamma_{\mathrm{N}})^d)
\end{eqnarray*}
\noindent
of $\Lambda$ in $(\lambda,\mu)\in L^\infty_+(\Omega)\times L^\infty_+(\Omega)$
fulfills
\begin{eqnarray}\label{Frechet}
\langle &\Lambda'(\lambda,\mu)(\hat{\lambda},\hat{\mu})g,h\rangle\\
&\quad=-\int_{\Omega} 2 \hat{\mu}\, \hat{\nabla}u^{g}_{(\lambda,\mu)} : \hat{\nabla}u^{h}_{(\lambda,\mu)} 
+ \hat{\lambda} \nabla \cdot u^{g}_{(\lambda,\mu)} \,\nabla\cdot  u^{h}_{(\lambda,\mu)}\,dx\nonumber,
\end{eqnarray}
where $\hat{\lambda},\hat{\mu}\in L^\infty(\Omega)$ and $u^{g}_{(\lambda,\mu)}, u^{h}_{(\lambda,\mu)}$ is the
solution of (\ref{direct_1}) for the boundary loads $g,h\in L^2(\Gamma_{\mathrm{N}})^d$, respectively.
\end{lemma}

\bproof
Clearly, it holds that for all $\lambda,\mu\in L^\infty_+(\Omega)$,
$\hat{\lambda},\hat{\mu}\in L^\infty(\Omega)$ fixed,
 $\Lambda^\prime(\lambda,\mu)(\hat{\lambda},\hat{\mu}):
 L^2(\Gamma_{\mathrm{N}})^d \to L^2(\Gamma_{\mathrm{N}})^d$ is linear and continuous and
equation (\ref{Frechet}) defines a continuous linear mapping
  $\Lambda^\prime(\lambda,\mu): (\hat{\lambda},\hat{\mu})\mapsto \Lambda^\prime(\lambda,\mu)(\hat{\lambda},\hat{\mu})$
 from $L^\infty(\Omega)\times L^\infty(\Omega)$ to $\mathcal{L}(L^2(\Gamma_{\mathrm{N}})^d)$
as well.
Further, it is obvious that $\Lambda^\prime(\lambda,\mu)(\hat{\lambda},\hat{\mu})$ is self-adjoint.
\\
Based on this, we will prove that $\Lambda$ is Fr\'echet differentiable and that $\Lambda^\prime$ is its Fr\'echet derivative, i.e.,
\begin{eqnarray*}
\fl \lim_{\Vert \tilde{h}\Vert_{L^\infty(\Omega)^2} \to 0}\frac{1}{\Vert \tilde{h}\Vert_{L^\infty(\Omega)^2} }
 \Vert\Lambda(\lambda+h^\lambda,\mu+h^\mu)-\Lambda(\lambda,\mu)-\Lambda^\prime(\lambda,\mu)(h^\lambda,h^\mu)\Vert_{\mathcal{L}(L^2(\Gamma_{\mathrm{N}})^d)}
 =0
\end{eqnarray*}
\noindent
for $\tilde{h}=(h^\lambda,h^\mu)^T$. Note that the space $L^\infty(\Omega)^2$ is equipped with the norm
\begin{eqnarray*}
\Vert \tilde{h}\Vert_{L^\infty(\Omega)^2}:=\max\left(\Vert h^\lambda\Vert_{L^\infty(\Omega)},\Vert h^\mu \Vert_{L^\infty(\Omega)}\right).
\end{eqnarray*}
\noindent
We take a look at 
\begin{eqnarray*}
&\fl\lim_{\Vert \tilde{h}\Vert_{L^\infty(\Omega)^2}  \to 0}\frac{1}{\Vert \tilde{h}\Vert_{L^\infty(\Omega)^2} }\Vert\Lambda(\lambda+h^\lambda,\mu+h^\mu)-\Lambda(\lambda,\mu)-\Lambda^\prime(\lambda,\mu)(h^\lambda,h^\mu)\Vert_{\mathcal{L}(L^2(\Gamma_{\mathrm{N}})^d)}\\
&\fl =\lim_{\Vert \tilde{h}\Vert_{L^\infty(\Omega)^2}  \to 0}\frac{1}{\Vert \tilde{h}\Vert_{L^\infty(\Omega)^2} }\sup_{\Vert g\Vert_{L^2(\Gamma_{\mathrm{N}})^d}=1}
\Big\vert\underbrace{\Big\langle g,\left(\Lambda(\lambda+h^\lambda,\mu+h^\mu)-\Lambda(\lambda,\mu)-\Lambda^\prime(\lambda,\mu)(h^\lambda,h^\mu)\right)g\Big\rangle}_{:=G}\Big\vert.
\end{eqnarray*}
\noindent
Next, we estimate $G$ via the monotonicity properties of the Neumann-to-Dirichlet operator $\Lambda(\lambda,\mu)$
(Lemma \ref{mono} and Lemma \ref{mon_est_2}):
\begin{eqnarray*}
G \geq& \int_{\Omega}\Bigg(2 (\mu-\mu-h^\mu)\hat{\nabla}u_{(\lambda,\mu)}^g : \hat{\nabla}u_{(\lambda,\mu)}^g
 +(\lambda-\lambda-h^\lambda)\nabla\cdot u_{(\lambda,\mu)}^g  \nabla\cdot u_{(\lambda,\mu)}^g\\
 &+ 2h^\mu \hat{\nabla}u_{(\lambda,\mu)}^g : \hat{\nabla}u_{(\lambda,\mu)}^g +
 h^\lambda \nabla\cdot u_{(\lambda,\mu)}^g : \nabla\cdot u_{(\lambda,\mu)}^g\Bigg)\,dx\\
\quad =&0
\end{eqnarray*}
\noindent
and
\begin{eqnarray*}
G \leq&\int_{\Omega}- 2\frac{\mu h^\mu}{\mu+h^\mu}\hat{\nabla}u_{(\lambda,\mu)}^g : \hat{\nabla}u_{(\lambda,\mu)}^g
 -\frac{\lambda h^\lambda}{\lambda+h^\lambda}\nabla\cdot u_{(\lambda,\mu)}^g  \nabla\cdot u_{(\lambda,\mu)}^g\\
 &+ 2h^\mu \hat{\nabla}u_{(\lambda,\mu)}^g : \hat{\nabla}u_{(\lambda,\mu)}^g +
 h^\lambda \nabla\cdot u_{(\lambda,\mu)}^g : \nabla\cdot u_{(\lambda,\mu)}^g\Bigg)\,dx\\
 \quad=&\int_{\Omega}\Bigg(2\frac{(h^\mu)^2}{\mu+h^\mu}\hat{\nabla}u_{(\lambda,\mu)}^g : \hat{\nabla}u_{(\lambda,\mu)}^g+
 \frac{(h^\lambda)^2}{\lambda+h^\lambda}\nabla\cdot u_{(\lambda,\mu)}^g  \nabla\cdot u_{(\lambda,\mu)}^g\Bigg)\,dx.
\end{eqnarray*}
Thus, we are led to
\begin{eqnarray*}
 &\fl \frac{1}{\Vert \tilde{h}\Vert_{L^\infty(\Omega)^2}}\int_{\Omega}\Bigg(2\frac{(h^\mu)^2}{\mu+h^\mu}\hat{\nabla}u_{(\lambda,\mu)}^g : \hat{\nabla}u_{(\lambda,\mu)}^g+
 \frac{(h^\lambda)^2}{\lambda+h^\lambda}\nabla\cdot u_{(\lambda,\mu)}^g  \nabla\cdot u_{(\lambda,\mu)}^g\Bigg)\,dx\\
&\fl \leq \frac{1}{\Vert \tilde{h}\Vert_{L^\infty(\Omega)^2}}
\Bigg(\Vert h^\mu\Vert_{L^\infty(\Omega)}\int_{\Omega}2\frac{\vert h^\mu\vert}{ \mu+h^\mu}\hat{\nabla}u_{(\lambda,\mu)}^g : \hat{\nabla}u_{(\lambda,\mu)}^g\,dx\\
&+\Vert h^\lambda\Vert_{L^\infty(\Omega)}\int_{\Omega}\frac{\vert h^\lambda\vert}{ \lambda+h^\lambda}\nabla\cdot u_{(\lambda,\mu)}^g  \nabla\cdot u_{(\lambda,\mu)}^g\,dx\Bigg)\\
&\fl\leq \frac{1}{\Vert \tilde{h}\Vert_{L^\infty(\Omega)^2}}
\underbrace{\max\left(\Vert h^\lambda\Vert_{L^\infty(\Omega)},\Vert h^\mu\Vert_{L^\infty(\Omega)}\right)}_{=\Vert\tilde{h}\Vert_{L^\infty(\Omega)^2}}
\Bigg(\int_{\Omega}2\frac{\vert h^\mu\vert}{ \mu+h^\mu}\hat{\nabla}u_{(\lambda,\mu)}^g : \hat{\nabla}u_{(\lambda,\mu)}^g\,dx\\
&+\int_{\Omega}\frac{\vert h^\lambda\vert}{ \lambda+h^\lambda}\nabla\cdot u_{(\lambda,\mu)}^g  \nabla\cdot u_{(\lambda,\mu)}^g\,dx\Bigg)\\
 &\fl=\int_{\Omega}2\frac{\vert h^\mu\vert}{ \mu+h^\mu}\hat{\nabla}u_{(\lambda,\mu)}^g : \hat{\nabla}u_{(\lambda,\mu)}^g\,dx+
 \int_{\Omega}\frac{\vert h^\lambda\vert}{\lambda+h^\lambda}\nabla\cdot u_{(\lambda,\mu)}^g  \nabla\cdot u_{(\lambda,\mu)}^g\,dx.
 \end{eqnarray*}
\noindent
If $\Vert \tilde{h}\Vert_{L^\infty(\Omega)^2}\to 0$, then $h^\lambda,h^\mu\to 0$ and we finally end up with
\begin{eqnarray*}
&\fl\lim_{\Vert \tilde{h}\Vert_{L^\infty(\Omega)^2}  \to 0}\frac{1}{\Vert \tilde{h}\Vert_{L^\infty(\Omega)^2} }
\Vert\Lambda(\lambda+h^\lambda,\mu+h^\mu)-\Lambda(\lambda,\mu)-\Lambda^\prime(\lambda,\mu)(h^\lambda,h^\mu)\Vert_{\mathcal{L}(L^2(\Gamma_{\mathrm{N}})^d)}
\to 0,
\end{eqnarray*}
\noindent
which proves that $(\hat{\lambda},\hat{\mu})\mapsto \Lambda^\prime(\lambda,\mu)(\hat{\lambda},\hat{\mu})$ is the Fr\'echet derivative 
of $\Lambda$ in $(\lambda,\mu)$.
\eproof

 \noindent
 \\
Finally, we want to state the following corollaries: 

\begin{corollary}\label{compact}
For $\lambda,\mu\in L^\infty_+(\Omega)$ and $\hat{\lambda},\hat{\mu}\in L^\infty(\Omega)$, 
the Fr\'echet derivative $\Lambda^\prime(\lambda,\mu)(\lambda,\mu)$ is compact and self-adjoint, and
for all $g\in L^2(\Gamma_{\mathrm{N}})^d$
\begin{eqnarray*}
\Lambda^\prime(\lambda,\mu)(\hat{\lambda},\hat{\mu})=v^g\vert_{\Gamma_{\mathrm{N}}},
\end{eqnarray*}
\noindent
where $v^g\in\mathcal{V}$ solves
\begin{eqnarray}\label{var_form_vg}
&\int_{\Omega}\left( 2\mu \hat{\nabla} v^g:\hat{\nabla}w+\lambda\nabla\cdot v^g \nabla\cdot w\right)\,dx\\
&=-\int_{\Omega}\left(2\hat{\mu}\hat{\nabla}u^g_{(\lambda,\mu)}:\hat{\nabla}w+\hat{\lambda} \nabla\cdot u^g_{(\lambda,\mu) \nabla\cdot w}\right)\,dx \nonumber
\end{eqnarray}
\noindent
for all $w\in\mathcal{V}$ and $u^g_{(\lambda,\mu)}\in\mathcal{V}$ solves (\ref{direct_1}).
\\
\\
Note that for $\mathrm{supp}(\lambda)\subseteq\Omega$, $\mathrm{supp}(\mu)\subseteq \Omega$, this is equivalent to
\begin{equation*}\hspace{-1.5cm}\left\{ \begin{array}{rcll}
\nabla\cdot\left(\lambda (\nabla\cdot v^g)I + 2\mu \hat{\nabla} v^g \right)&=&
-\nabla\cdot\left(\hat{\lambda} \left(\nabla\cdot u_{(\lambda,\mu)}^g\right)I + 2\hat{\mu} \hat{\nabla} u_{(\lambda,\mu)}^g \right)\,\, \textnormal{\it in}\,\,\, \Omega,\\
\left(\lambda (\nabla\cdot v^g)I + 2\mu \hat{\nabla} v^g \right) n&=&0\,\,\,\textnormal{\it on}\,\,\, \Gamma_{\mathrm{N}},\\
v^g&=&0\,\,\, \textnormal{\it on}\,\,\, \Gamma_{\mathrm{D}}.
\end{array}\right.
\end{equation*}
\end{corollary}

\bproof
The unique solvability of (\ref{var_form_vg}) follows form the Lax-Milgram theorem. Using (\ref{var_form_vg})
together with the weak formulation (\ref{var-direct_1}) shows that
\begin{eqnarray*}
 \langle v^g\vert_{\Gamma_{\mathrm{N}}},h\rangle
 &=\int_{\Omega}\left(2\mu \hat{\nabla}v^g:\hat{\nabla}u^h_{(\lambda,\mu)}
 +\lambda \nabla\cdot v^g\nabla\cdot u^h_{(\lambda,\mu)}\right)\,dx\\
 &=-\int_{\Omega}\left(2\hat{\mu} \hat{\nabla}u^g_{(\lambda,\mu)}:u^h_{(\lambda,\mu)}
 +\hat{\lambda}\nabla\cdot u^g_{(\lambda,\mu)} \nabla\cdot u^h_{(\lambda,\mu)}\right)\,dx\\
 &=\langle \Lambda^\prime(\lambda,\mu)(\hat{\lambda},\hat{\mu})g,h\rangle
\end{eqnarray*}
\noindent
for all $h\in L^2(\Gamma_{\mathrm{N}})^d$.
\\
\\
Self-adjointness follows from (\ref{Frechet}) and compactness follows from the compactness of the trace operator.
\eproof

\begin{corollary}
For  $\hat{\lambda}_0, \hat{\lambda}_1, \hat{\mu}_0, \hat{\mu}_1 \in L^\infty(\Omega)$, we have that for
\begin{eqnarray*}
\hat{\lambda}_0\leq \hat{\lambda}_1 \,\,\textnormal{\it and }\,\,  \hat{\mu}_0 \leq \hat{\mu}_1
\,\,
\textnormal{\it it follows that}\,\,
\Lambda'(\lambda,\mu)(\hat{\lambda}_0,\hat{\mu}_0)\geq \Lambda'(\lambda,\mu)(\hat{\lambda}_1,\hat{\mu}_1).
\end{eqnarray*}
\end{corollary}

\bproof
Let $\hat{\lambda}_0\leq \hat{\lambda}_1$ and  $\hat{\mu}_0 \leq \hat{\mu}_1$. Then
\begin{eqnarray*}
 \langle \Lambda^\prime&(\lambda,\mu)(\hat{\lambda}_0,\hat{\mu}_0)g,g\rangle
 -\langle \Lambda^\prime(\lambda,\mu)(\hat{\lambda}_1,\hat{\mu}_1)g,g\rangle\\
 =&-\int_{\Omega}2\hat{\mu}_0 \hat{\nabla}u_{(\lambda,\mu)}^g : \hat{\nabla}u_{(\lambda,\mu)}^g
 +\hat{\lambda}_0 \nabla\cdot u_{(\lambda,\mu)}^g  \nabla\cdot u_{(\lambda,\mu)}^g\,dx\\
 &+\int_{\Omega}2\hat{\mu}_1 \hat{\nabla}u_{(\lambda,\mu)}^g : \hat{\nabla}u_{(\lambda,\mu)}^g
 +\hat{\lambda}_1 \nabla\cdot u_{(\lambda,\mu)}^g  \nabla\cdot u_{(\lambda,\mu)}^g\,dx\\
 \geq& 0.
\end{eqnarray*}
\eproof
}
\noindent
\\
Based on these results, we can now introduce the standard monotonicity method.

\subsection{Standard monotonicity method}
Our aim is to prove that the opposite direction of Corollary \ref{monotonicity} holds true in order to formulate the so-called standard monotonicity test 
(Corollary \ref{col_stan_mon_1} and \ref{col_stan_mon_2}).
\\
\\
{\color{black}
Therefore, we consider the case where $\Omega$ contains inclusions in which the Lam\'e parameters $(\lambda,\mu)$ differ from otherwise 
known constant background Lam\'e  parameters $(\lambda_0,\mu_0)$. 
{\color{black}
We will derive a method to reconstruct the unknown inclusion that is based on comparing the measured Neumann-to-Dirichlet 
operator $\Lambda(\lambda,\mu)$ with that of a test inclusion.
}
}
\\
For the precise formulation, we will now introduce the concept of the inner and the outer support of a measurable function in a similar way as in 
\cite{harrach2013monotonicity}.
\begin{definition}
A relatively open set  $U\subseteq \overline{\Omega}$ is called connected to $\partial\Omega$ if  $U \cap\Omega$
is connected and $U\cap \partial \Omega \neq \emptyset$.
\end{definition}

\noindent

\begin{definition}
 For a measurable function  {\color{black} $\varphi=(\varphi_1,\varphi_2): \Omega \rightarrow \mathbb{R}^2$}, we define
 \begin{itemize}
 \item[a)] the support $\mathrm{supp}(\varphi)$ as the complement (in $\overline{\Omega}$) of the union of those relatively
 open $U\subseteq\overline{\Omega}$, for which  $\varphi|_{U}\equiv 0$,
 \item[b)] 
 {\color{black}
 the inner support $\mathrm{inn}\,\mathrm{supp}(\varphi_i)$ as the union of those open sets $U\subseteq \Omega$, for which  
 $\mathrm{ess}\,\inf_{x\in U}\vert\varphi_i(x)\vert >0$, $i=1,2$, and 
 $\mathrm{inn}\,\mathrm{supp}(\varphi):=\mathrm{inn}\,\mathrm{supp}(\varphi_1)\cup\mathrm{inn}\,\mathrm{supp}(\varphi_2)$,
 }
 \item[c)] the outer support $\mathrm{out}_{\partial\Omega}\, \mathrm{supp}(\varphi)$ as the complement (in $\overline{\Omega}$) of the union of those 
 relatively open $U\subseteq \overline{\Omega}$ that are connected to $\partial\Omega$ and for which 
 $\varphi|_{U}\equiv 0$.
 \end{itemize}
\end{definition}
{\color{black}
We now consider 
$\varphi=(\lambda-\lambda_0, \mu-\mu_0)^T$.
Our goal is to recover $\mathrm{supp}(\varphi)$
from the knowledge of the Neumann-to-Dirichlet operator $\Lambda(\lambda,\mu)$.
We show that this can be done up to the difference between the inner and outer support.
\\
Let us consider the setting as depicted in Figure \ref{plot_support}. By proving the opposite direction of Corollary \ref{monotonicity},
we show that $\mathrm{supp}(\varphi)$ can be reconstructed by monotonicity
tests, which simply compare $\Lambda(\lambda,\mu)$ (in the sense of quadratic forms) to  the Neumann-to-Dirichlet 
operators of test  parameters.
To be more precise, the support of $\varphi$
can be reconstructed under the assumption that $\mathrm{supp}(\varphi)\subset \Omega$ has a connected complement, in which case we have
$\mathrm{supp}(\varphi)=\mathrm{out}_{\partial\Omega}\, \mathrm{supp}(\varphi)$ (c.f. \cite{harrach2013monotonicity}).
Otherwise, what we can reconstruct is the support of $\varphi$ together
with all holes that have no connection to the boundary  $\partial\Omega$ (c.f. \cite{Factorization_Harrach_2013}). 
 \begin{figure}[h]
 \begin{center}
 \includegraphics[width=0.6\textwidth]{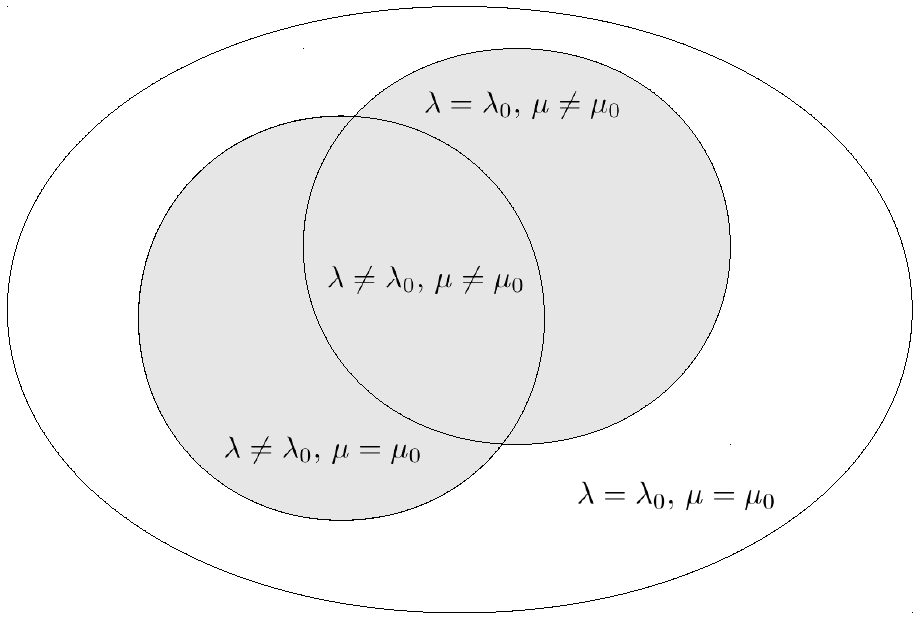}
\caption{Exemplary $\mathrm{out}_{\partial\Omega}\, \mathrm{supp}(\varphi)$ for an inclusion without holes in gray.}\label{plot_support}
\end{center}
 \end{figure}
 \noindent
 \\
In the following, we define $\alpha$ and $\beta$ as the contrasts and $\chi_{\mathcal{D}}$ as well as
$\chi_\mathcal{B}$ as the characteristic function w.r.t. the inclusion 
$\mathcal{D}:=\mathrm{out}_{\partial\Omega}\,\mathrm{supp}(\varphi)$ and 
the so-called test inclusion $\mathcal{B}$, respectively. 

%
%

\noindent
\begin{theorem}\label{theor_inclusion_1}
Let $\lambda_0,\mu_0>0$ be constant and let $(\lambda,\mu)\in {\color{black}L^\infty_+(\Omega)\times L^\infty_+(\Omega)}$, 
$\lambda\geq \lambda_0$, $\mu\geq \mu_0$.
For every open set
$\mathcal{B}$ (e.g. ball or cube) and every $\alpha,\beta\geq 0$, $\alpha+\beta> 0$,
\begin{eqnarray*}
\alpha \chi_{\mathcal{B}}\leq \lambda -\lambda_0,\quad\beta \chi_{\mathcal{B}}\leq \mu -\mu_0 
\end{eqnarray*}
\noindent
implies
\begin{eqnarray*}
\Lambda(\lambda_0+\alpha \chi_{\mathcal{B}}, \mu_0+\beta \chi_{\mathcal{B}})\geq \Lambda(\lambda,\mu)
\end{eqnarray*}
\noindent
and
\begin{eqnarray*}
\mathcal{B}\nsubseteq \mathcal{D} \,\,\textnormal{\it implies}\,\,
\Lambda(\lambda_0+\alpha \chi_{\mathcal{B}}, \mu_0+\beta \chi_{\mathcal{B}})\ngeq \Lambda(\lambda,\mu).
\end{eqnarray*}
\noindent
Hence, the set
\begin{eqnarray*}
R:=\bigcup\limits_{\alpha,\beta\geq 0, \alpha+\beta> 0}\left\{ \mathcal{B}\subseteq\Omega:
\Lambda(\lambda_0 +\alpha \chi_{\mathcal{B}}, \mu_0 +\beta\chi_{\mathcal{B}})\geq \Lambda(\lambda,\mu)\right\},
\end{eqnarray*}
\noindent
fulfills
\begin{eqnarray*}
\mathrm{inn }\,\mathrm{supp}((\lambda-\lambda_0,\mu - \mu_0)^T)
\subseteq R \subseteq 
\mathcal{D}.
\end{eqnarray*}
\end{theorem}

\bproof
Let $\lambda_0,\mu_0>0$ be constant and let $(\lambda,\mu)\in {\color{black}L^\infty_+(\Omega)\times L^\infty_+(\Omega)}$, 
$\lambda\geq \lambda_0$, $\mu\geq \mu_0$.
 Let $\mathcal{B}$ be an open set and $\alpha,\beta\geq 0$, $\alpha+\beta>0$ and
\begin{eqnarray*}
\alpha \chi_{\mathcal{B}}\leq \lambda-\lambda_0,\quad \beta \chi_{\mathcal{B}}\leq \mu -\mu_0.
\end{eqnarray*}
\noindent
We immediately obtain from Corollary \ref{monotonicity} that
\begin{eqnarray*}
\Lambda(\lambda_0+\alpha \chi_{\mathcal{B}}, \mu_0+\beta \chi_{\mathcal{B}})\geq \Lambda(\lambda,\mu).
\end{eqnarray*}
\noindent
It remains to show that 
\begin{eqnarray}\label{B_not_in_D}
\mathcal{B} \nsubseteq\mathcal{D} \quad \mathrm{ implies } \quad
\Lambda(\lambda_0+\alpha \chi_{\mathcal{B}}, \mu_0+\beta \chi_{\mathcal{B}})\ngeq \Lambda(\lambda,\mu).
\end{eqnarray}
\noindent
Let $\mathcal{B}\subseteq \Omega$ with $\mathcal{B}\not\subseteq \mathcal{D}$. 
Then there exists a non-empty open set $O$ with $\overline{O}\subset \mathcal{B}\setminus\mathcal{D}$ and 
we introduce a function $\phi\in C^\infty(\Omega)$ which fulfills  $0\leq\phi\leq \chi_{\mathcal{B}\setminus\mathcal{D}}$ 
and $\phi\geq \frac{1}{2}$ 
on $O$.
\\
We now apply Theorem 3.3 from \cite{Eberle_Monotonicity} and obtain a sequence $(g_n)_{n\in\mathbb{N}}\subset L^2( \Gamma_{\mathrm{N}})^d$ 
so that the solutions 
\mbox{$(u^{(g_n)})_{n\in\mathbb{N}}\subset \mathcal{V}$}  of
\begin{equation}\hspace{-1cm}\left\{ \begin{array}{rcll}\label{theo_3_a}
\nabla\cdot\left((\lambda_0+\alpha \phi) (\nabla\cdot u^{(g_n)})I + 2(\mu_0+\beta \phi) \hat{\nabla} u^{(g_n)} \right)&=&0 &\mathrm{in}\,\,\Omega,\\
\left((\lambda_0+\alpha \phi) (\nabla\cdot u^{(g_n)})I + 2(\mu_0+\beta \phi) \hat{\nabla} u^{(g_n)} \right) n&=&g_n&\mathrm{on}\,\, \Gamma_{\mathrm{N}},\\
u^{(g_n)}&=&0 &\mathrm{on}\,\, \Gamma_{\mathrm{D}},
\end{array}\right.
\end{equation}
fulfill
\begin{eqnarray}
&\lim_{n\to \infty}\int_{O} \hat{\nabla}u^{(g_n)}:\hat{\nabla}u^{(g_n)}\, dx=\infty,\label{loc_pot_1_a}\\
&\lim_{n\to \infty}\int_{\mathcal{D}} \hat{\nabla}u^{(g_n)}:\hat{\nabla}u^{(g_n)}\, dx=0,\label{loc_pot_2_a}\\
&\lim_{n\to \infty}\int_{O} \nabla\cdot u^{(g_n)}\, \nabla \cdot u^{(g_n)}\, dx=\infty,\label{loc_pot_3_a}\\
&\lim_{n\to \infty}\int_{\mathcal{D}}  \nabla\cdot u^{(g_n)}\, \nabla \cdot u^{(g_n)}\, dx=0\label{loc_pot_4_a}.
\end{eqnarray}

\noindent
From Corollary \ref{monotonicity}, Lemma \ref{mono} and Equation (\ref{loc_pot_1_a})-(\ref{loc_pot_4_a}) it follows 
\begin{eqnarray*}
&\langle g_n, \Lambda(\lambda_0+\alpha \chi_{\mathcal{B}},\mu_0+\beta \chi_{\mathcal{B}}),g_n\rangle
-\langle g_n, \Lambda(\lambda,\mu)g_n\rangle\\
&\quad\leq \langle g_n, \Lambda(\lambda_0+\alpha\phi,\mu_0+\beta\phi)g_n\rangle-\langle g_n,\Lambda(\lambda,\mu)g_n\rangle\\
&\quad\leq \int_\Omega 2(\mu-\mu_0-\beta\phi)\hat{\nabla}u^{(g_n)}:\hat{\nabla}u^{(g_n)}\\
&\qquad+(\lambda-\lambda_0-\alpha\phi)\nabla\cdot u^{(g_n)} \nabla\cdot u^{(g_n)}\,dx\\
&\quad=\int_{\mathcal{D}} \left(2(\mu-\mu_0)\hat{\nabla}u^{(g_n)}:\hat{\nabla}u^{(g_n)} 
+(\lambda-\lambda_0)\nabla\cdot u^{(g_n)} \nabla\cdot u^{(g_n)}\right)\,dx\\
&\qquad-\int_{\Omega}\left(2\beta\phi\hat{\nabla}u^{(g_n)}:\hat{\nabla}u^{(g_n)}+\alpha\phi\nabla \cdot u^{(g_n)}\nabla\cdot u^{(g_n)}\right)\,dx\\
&\quad\leq\int_{\mathcal{D}}\left( 2 (\mu-\mu_0) \hat{\nabla}u^{(g_n)} : \hat{\nabla}u^{(g_n)}+(\lambda-\lambda_0)\nabla\cdot u^{(g_n)} \nabla\cdot u^{(g_n)}\right)\,dx\\
&\qquad -\frac{1}{2}\int_{O} \left(2 \beta \hat{\nabla}u^{(g_n)} : \hat{\nabla}u^{(g_n)}+\alpha  \nabla\cdot u^{(g_n)} \nabla\cdot u^{(g_n)}\right)\,dx\\
&\quad\to -\infty,
\end{eqnarray*}
\noindent
and hence
\begin{eqnarray*}
\Lambda(\lambda_0+\alpha \chi_{\mathcal{B}}, \mu_0+\beta \chi_{\mathcal{B}})\ngeq \Lambda(\lambda,\mu).
\end{eqnarray*}
\eproof
}
\\
{\color{black}
In addition, we state the theorem for the case $\lambda\leq \lambda_0$, $\mu\leq \mu_0$.

\begin{theorem}\label{theor_inclusion_add}
Let $\lambda_0,\mu_0>0$ be constant and let $(\lambda,\mu)\in {\color{black}L^\infty_+(\Omega)\times L^\infty_+(\Omega)}$, 
$\lambda\leq \lambda_0$, $\mu\leq \mu_0$.
For every open set
$\mathcal{B}$ (e.g. ball or cube) and every $\alpha,\beta\geq 0$, $\alpha+\beta> 0$,
\begin{eqnarray*}
\alpha \chi_{\mathcal{B}}\leq \lambda_0 -\lambda,\quad\beta \chi_{\mathcal{B}}\leq \mu_0 -\mu 
\end{eqnarray*}
\noindent
implies
\begin{eqnarray*}
\Lambda(\lambda_0-\alpha \chi_{\mathcal{B}}, \mu_0-\beta \chi_{\mathcal{B}})\leq \Lambda(\lambda,\mu)
\end{eqnarray*}
\noindent
and
\begin{eqnarray*}
\mathcal{B}\nsubseteq \mathcal{D}
\,\, \textnormal{\it implies}\,\,
\Lambda(\lambda_0-\alpha \chi_{\mathcal{B}}, \mu_0-\beta \chi_{\mathcal{B}})\nleq \Lambda(\lambda,\mu).
\end{eqnarray*}
\noindent
Hence, the set
\begin{eqnarray*}
R:=\bigcup\limits_{\alpha,\beta\geq 0, \alpha+\beta> 0}\left\{ \mathcal{B}\subseteq\Omega: \Lambda(\lambda_0 -\alpha \chi_{\mathcal{B}}, \mu_0 -\beta\chi_{\mathcal{B}})\leq \Lambda(\lambda,\mu)\right\},
\end{eqnarray*}
\noindent
fulfills
\begin{eqnarray*}
\mathrm{inn }\,\mathrm{supp}((\lambda_0 -\lambda,\mu_0-\mu)^T)
\subseteq R \subseteq 
\mathcal{D}.
\end{eqnarray*}
\end{theorem}

\bproof
Let $\lambda_0,\mu_0>0$ be constant and let $(\lambda,\mu)\in {\color{black}L^\infty_+(\Omega)\times L^\infty_+(\Omega)}$, 
$\lambda\leq \lambda_0$, $\mu\leq \mu_0$.
 Let $\mathcal{B}$ be an open set, $\alpha,\beta\geq 0$, $\alpha+\beta>0$ and
\begin{eqnarray*}
\alpha \chi_{\mathcal{B}}\leq \lambda_0-\lambda,\quad \beta \chi_{\mathcal{B}}\leq \mu_0 -\mu.
\end{eqnarray*}
Hence, we get with Corollary \ref{monotonicity} that
\begin{eqnarray*}
\Lambda(\lambda_0-\alpha \chi_{\mathcal{B}}, \mu_0-\beta \chi_{\mathcal{B}})\leq \Lambda(\lambda,\mu).
\end{eqnarray*}
\noindent
It remains to show that 
\begin{eqnarray*}
\mathcal{B} \nsubseteq\mathcal{D} \quad \mathrm{ implies } \quad
\Lambda(\lambda_0-\alpha \chi_{\mathcal{B}}, \mu_0-\beta \chi_{\mathcal{B}})\nleq \Lambda(\lambda,\mu).
\end{eqnarray*}
Let $\mathcal{B}\nsubseteq \mathcal{D}$.
Similar as in the proof of Theorem \ref{theor_inclusion_1}, we
introduce a function $\phi\in C^\infty(\Omega)$ which fulfills  $0\leq\phi\leq \chi_{\mathcal{B}\setminus\mathcal{D}}$ and
$\phi\geq \frac{1}{2}$ on a non-empty open set $O$ with $\overline{O}\subset \mathcal{B}\setminus\mathcal{D}$.
We now apply Theorem 3.3 from \cite{Eberle_Monotonicity} and obtain a sequence $(g_n)_{n\in\mathbb{N}}\subset L^2( \Gamma_{\mathrm{N}})^d$ so that the solutions 
\mbox{$(u^{(g_n)})_{n\in\mathbb{N}}\subset \mathcal{V}$}  of
\begin{equation}\hspace{-1cm}\left\{ \begin{array}{rcll}\label{theo_3_a_new}
\nabla\cdot\left((\lambda_0-\alpha \phi ) (\nabla\cdot u^{(g_n)})I + 2(\mu_0-\beta \phi ) \hat{\nabla} u^{(g_n)} \right)&=&0 &\mathrm{in}\,\,\Omega,\\
\left((\lambda_0-\alpha \phi ) (\nabla\cdot u^{(g_n)})I + 2(\mu_0-\beta \phi ) \hat{\nabla} u^{(g_n)} \right) n&=&g_n&\mathrm{on}\,\, \Gamma_{\mathrm{N}},\\
u^{(g_n)}&=&0 &\mathrm{on}\,\, \Gamma_{\mathrm{D}},
\end{array}\right.
\end{equation}
fulfill
\begin{eqnarray}
&\lim_{n\to \infty}\int_{O} \hat{\nabla}u^{(g_n)}:\hat{\nabla}u^{(g_n)}\, dx=\infty,\label{loc_pot_1}\\
&\lim_{n\to \infty}\int_{\mathcal{D}} \hat{\nabla}u^{(g_n)}:\hat{\nabla}u^{(g_n)}\, dx=0,\label{loc_pot_2}\\
&\lim_{n\to \infty}\int_{O} \nabla\cdot u^{(g_n)}\, \nabla \cdot u^{(g_n)}\, dx=\infty,\label{loc_pot_3}\\
&\lim_{n\to \infty}\int_{\mathcal{D}}  \nabla\cdot u^{(g_n)}\, \nabla \cdot u^{(g_n)}\, dx=0\label{loc_pot_4}.
\end{eqnarray}

\noindent
From Corollary \ref{mono}, Lemma \ref{mon_est_2} and Equation (\ref{loc_pot_1})-(\ref{loc_pot_4}) it follows 
\begin{eqnarray*}
&\langle g_n, \Lambda(\lambda_0-\alpha \chi_{\mathcal{B}},\mu_0-\beta \chi_{\mathcal{B}}),g_n\rangle
-\langle g_n, \Lambda(\lambda,\mu)g_n\rangle\\
&\quad\geq \langle g_n, \Lambda(\lambda_0-\alpha\phi ,\mu_0-\beta\phi )g_n\rangle-\langle g_n,\Lambda(\lambda,\mu)g_n\rangle\\
&\quad\geq \int_{\Omega} \left(2\frac{\mu_0-\beta\phi }{\mu}(\mu-\mu_0+\beta\phi )\hat{\nabla}u^{(g_n)}:\hat{\nabla}u^{(g_n)}\right)\,dx\\
&\qquad +\int_{\Omega}\left(\frac{\lambda_0-\alpha\phi }{\lambda}(\lambda-\lambda_0+\alpha\phi ) \nabla\cdot u^{(g_n)} \nabla\cdot u^{(g_n)}\right)\,dx\\
&\quad\geq \int_{\mathcal{D}}\left( 2\frac{\mu_0}{\mu}(\mu-\mu_0)\hat{\nabla}u^{(g_n)}:\hat{\nabla}u^{(g_n)}
 +\frac{\lambda_0}{\lambda}(\lambda-\lambda_0) \nabla \cdot u^{(g_n)} \nabla \cdot u^{(g_n)}\right)\,dx\\
&\qquad +\frac{1}{2}\int_{O} \left(2\frac{\mu_0-\beta }{\mu_0}\beta \hat{\nabla}u^{(g_n)}:\hat{\nabla}u^{(g_n)}
+\frac{\lambda_0-\alpha }{\lambda_0}\alpha\nabla \cdot u^{(g_n)} \nabla \cdot u^{(g_n)}\right)\,dx\\
&\quad \geq \infty
\end{eqnarray*}
\noindent
and hence
\begin{eqnarray*}
\Lambda(\lambda_0-\alpha \chi_{\mathcal{B}}, \mu_0-\beta \chi_{\mathcal{B}})\nleq \Lambda(\lambda,\mu).
\end{eqnarray*}
\eproof
}
\\
Based on this, we introduce the monotonicity tests. 
{\color{black}
In what follows, we reconstruct an inclusion $\mathcal{D}$ under the assumption, that
the complement of $\mathcal{D}\subset{\Omega}$ is connected to the boundary $\partial\Omega$.
In addition, we assume that the Lam\'e parameters $(\lambda_0,\mu_0)$ of the background as well 
as the Lam\'e parameters $(\lambda_1,\mu_1)$
inside the inclusion are constant.

}
\begin{corollary}{Standard monotonicity test: 1. version}\label{col_stan_mon_1}
\\
{\color{black}
Let $\lambda_0$, $\lambda_1$, $\mu_0$, $\mu_1\in\mathbb{R}^+$ with $\lambda_1>\lambda_0$, $\mu_1>\mu_0$  
and assume that 
\mbox{$(\lambda,\mu)=$} \mbox{$(\lambda_0+(\lambda_1-\lambda_0)\chi_\mathcal{D},\mu_0+(\mu_1-\mu_0)\chi_{\mathcal{D}})$}
with $\mathcal{D}=\mathrm{out}_{\partial\Omega}\,\mathrm{supp}((\lambda-\lambda_0,\mu-\mu_0)^T)$.}
Further on, let $\alpha,\beta\geq 0$, $\alpha+\beta>0$ with $\alpha \leq \lambda_1 -\lambda_0$, $\beta\leq \mu_1 -\mu_0$.
Then for every open set $\mathcal{B}\subseteq\Omega$
\begin{eqnarray}\label{stand_mono_test}
\hspace{-1.5cm}
\mathcal{B}\subseteq\mathcal{D}\quad\mathrm{if\,and\,only\,if}\quad\Lambda(\lambda_0+\alpha\chi_\mathcal{B},\mu_0+\beta\chi_\mathcal{B})\geq\Lambda(\lambda,\mu).
\end{eqnarray}
\end{corollary}
{\color{black}
\bproof
Let $\mathcal{B}\subseteq \mathcal{D}$ and $\alpha \leq \lambda_1 -\lambda_0$, $\beta\leq \mu_1 -\mu_0$. This means, that the conditions 
$\alpha \chi_\mathcal{B}\leq \lambda -\lambda_0$, $\beta \chi_\mathcal{B}\leq \mu -\mu_0 $ of Theorem \ref{theor_inclusion_1} are fulfilled
which immediately leads to $\Lambda(\lambda_0+\alpha \chi_\mathcal{B}, \mu_0+\beta \chi_\mathcal{B})\geq \Lambda(\lambda,\mu)$.
\\
For the opposite direction we assume that there exits a $\mathcal{B}\nsubseteq \mathcal{D}$, which fulfills $\Lambda(\lambda_0+\alpha \chi_\mathcal{B}, \mu_0+\beta \chi_\mathcal{B})\geq \Lambda(\lambda,\mu)$.
By applying the second part of Theorem \ref{theor_inclusion_1} we obtain
$\Lambda(\lambda_0+\alpha \chi_\mathcal{B}, \mu_0+\beta \chi_\mathcal{B})\ngeq \Lambda(\lambda,\mu)$ which contradicts 
$\Lambda(\lambda_0+\alpha \chi_\mathcal{B}, \mu_0+\beta \chi_\mathcal{B})\geq \Lambda(\lambda,\mu).$
Hence, $\Lambda(\lambda_0+\alpha \chi_\mathcal{B}, \mu_0+\beta \chi_\mathcal{B})\geq \Lambda(\lambda,\mu)$ implies that $\mathcal{B}\subseteq \mathcal{D}$.
\eproof
}
\noindent
\\
\\
Further on, we formulate the corresponding corollary for the case $\lambda\leq \lambda_0$ and $\mu\leq \mu_0$.
\begin{corollary}{Standard monotonicity test: 2. version}\label{col_stan_mon_2}
\\
{\color{black}
Let $\lambda_0$, $\lambda_1$, $\mu_0$, $\mu_1\in\mathbb{R}^+$ with $\lambda_1<\lambda_0$, $\mu_1<\mu_0$  
and assume that 
\mbox{$(\lambda,\mu)=$} \mbox{$(\lambda_0+(\lambda_1-\lambda_0)\chi_\mathcal{D},\mu_0+(\mu_1-\mu_0)\chi_{\mathcal{D}})$}
with $\mathcal{D}=\mathrm{out}_{\partial\Omega}\,\mathrm{supp}((\lambda-\lambda_0,\mu-\mu_0)^T)$.}
Further on, let $\alpha,\beta\geq 0$, $\alpha+\beta>0$ with $\alpha \leq \lambda_0 -\lambda_1$, $\beta\leq \mu_0 -\mu_1$.
Then for every open set $\mathcal{B}\subseteq\Omega$
\begin{eqnarray}\label{stand_mono_test}
\hspace{-1.5cm}
\mathcal{B}\subseteq\mathcal{D}\quad\mathrm{if\,and\,only\,if}\quad\Lambda(\lambda_0-\alpha\chi_\mathcal{B},\mu_0-\beta\chi_\mathcal{B})\leq\Lambda(\lambda,\mu).
\end{eqnarray}
\end{corollary}
{\color{black}
\bproof
The proof follows the lines of the proof of Corollary \ref{col_stan_mon_1} but we have to consider Theorem \ref{theor_inclusion_add}, where 
$\alpha \chi_\mathcal{B}\leq \lambda_0 -\lambda$, $\beta \chi_\mathcal{B}\leq \mu_0 -\mu$ 
implies
$\Lambda(\lambda_0-\alpha \chi_\mathcal{B}, \mu_0-\beta \chi_\mathcal{B})\leq \Lambda(\lambda,\mu)$.
\\
Further on, we use that $\mathcal{B}\nsubseteq \mathcal{D}$
implies $\Lambda(\lambda_0-\alpha \chi_\mathcal{B}, \mu_0-\beta \chi_\mathcal{B})\nleq \Lambda(\lambda,\mu)$.
\eproof
}
\noindent
\\
\\
Next, we apply Theorem \ref{theor_inclusion_1} to difference measurements 
\begin{eqnarray}
\Lambda_\mathcal{D}=\Lambda(\lambda_0,\mu_0)-\Lambda(\lambda,\mu)\label{Lambda_D}
\end{eqnarray}
\noindent which leads directly to the following lemma.

\begin{lemma}\label{mono_test_D}
Let $\Lambda_\mathcal{B}=\Lambda(\lambda_0,\mu_0)-\Lambda(\lambda_0+\alpha\chi_\mathcal{B},\mu_0+\beta \chi_\mathcal{B})$. Under the same assumptions on $\lambda$ and $\mu$
as in Theorem \ref{theor_inclusion_1}, we have for every open set $\mathcal{B}$
(e.g. ball or cube) and every $\alpha,\beta\geq 0$, $\alpha+\beta>0$,
\begin{eqnarray}\label{alpha_beta}
\alpha \chi_\mathcal{B}\leq \lambda -\lambda_0,\quad\beta \chi_\mathcal{B}\leq \mu -\mu_0 
\end{eqnarray}
\noindent
implies
\begin{eqnarray}\label{mon_diff_meas}
\Lambda_\mathcal{D}-\Lambda_\mathcal{B} \geq 0
\end{eqnarray}
\noindent
and
\begin{eqnarray*}
\mathcal{B}\nsubseteq\mathcal{D}
\,\,\textnormal{\it implies}\,\,
\Lambda_\mathcal{D}-\Lambda_\mathcal{B}\ngeq 0.
\end{eqnarray*}
\end{lemma}
\noindent
\\
{\color{black}
\bproof 
We take a look at the difference $\Lambda_\mathcal{D}-\Lambda_\mathcal{B}=\Lambda(\lambda_0+\alpha\chi_\mathcal{B},\mu_0+\beta\chi_\mathcal{B})-\Lambda(\lambda,\mu)$.
The assumption $\alpha\chi_\mathcal{B}\leq \lambda-\lambda_0$, $\beta\chi_\mathcal{B}\leq \mu-\mu_0$ leads via Theorem \ref{theor_inclusion_1} directly to the desired results.
\eproof
}
\noindent 
\\
\\
Next, we go over to the consideration of noisy difference measurements 
\begin{eqnarray}\label{Lambda_delta}
 \Lambda^\delta\approx\Lambda_\mathcal{D}
\end{eqnarray}
\noindent
with 
\begin{eqnarray}\label{Lambda_delta_norm}
\Vert{\Lambda^\delta - \Lambda_\mathcal{D}}\Vert<\delta,
\end{eqnarray}
\noindent
where $\delta>0$ and formulate a monotonicity test (c.f$.$ Corollary \ref{col_stan_mon_1}).
{\color{black}
By replacing $\Lambda^\delta-\Lambda_\mathcal{B}$ by its symmetric part, e.g. by $((\Lambda^\delta-\Lambda_\mathcal{B})+(\Lambda^\delta-\Lambda_\mathcal{B})^*)/2$, 
we can assume without loss of generality that $\Lambda^\delta-\Lambda_\mathcal{B}$ is self-adjoint.
}
{\color{black} 
Thus, the standard monotonicity test for noisy difference measurements reads as follows.
}
{\color{black}
\begin{corollary}{Standard monotonicity test for noisy difference measurements}\label{standard_test_noise}
\\
{\color{black}
Let $\lambda_0$, $\lambda_1$, $\mu_0$, $\mu_1\in\mathbb{R}^+$ with $\lambda_1>\lambda_0$, $\mu_1>\mu_0$  
and assume that 
\mbox{$(\lambda,\mu)=$} \mbox{$(\lambda_0+(\lambda_1-\lambda_0)\chi_\mathcal{D},\mu_0+(\mu_1-\mu_0)\chi_{\mathcal{D}})$}
with $\mathcal{D}=\mathrm{out}_{\partial\Omega}\,\mathrm{supp}((\lambda-\lambda_0,\mu-\mu_0)^T)$.}
Further on, let $\alpha,\beta\geq 0$, $\alpha+\beta>0$ with $\alpha \leq \lambda_1 -\lambda_0$, $\beta\leq \mu_1 -\mu_0$.
{\color{black}
Let $\Lambda^\delta$ be the Neumann-to-Dirichlet operator for noisy difference measurements with 
noise level $\delta>0$.
Then for every open set $\mathcal{B}\subseteq\Omega$ there exists a noise level $\delta_0>0$,
such that for all $0<\delta< \delta_0$, $\mathcal{B}$ is correctly detected as inside 
{\color{black} or not inside the inclusion $\mathcal{D}$ by the following monotonicity test
\begin{eqnarray}\label{test_standard_noise}
\mathcal{B}\subseteq\mathcal{D}\quad \textnormal{\it if and only if} \quad\Lambda^\delta-\Lambda_\mathcal{B} +\delta I\geq 0.
\end{eqnarray}
}
}
\end{corollary}
\bproof
Let $\mathcal{B}\not\subseteq \mathcal{D}$. $\Lambda_{\mathcal{D}}-\Lambda_{\mathcal{B}}$
is compact and self-adjoint (c.f. Corollary \ref{Lambda_compact}) and by Lemma \ref{mono_test_D} $\Lambda_{\mathcal{D}}-\Lambda_{\mathcal{B}}\not\geq 0$.
Hence, $\Lambda_{\mathcal{D}}-\Lambda_{\mathcal{B}}$ has negative eigenvalues. Let $\theta<0$ be the smallest eigenvalue with corresponding
eigenvector $g\in L^2(\Gamma_{\mathrm{N}})^d$, so that
\begin{eqnarray*}
 \langle (\Lambda^\delta-\Lambda_{\mathcal{B}}+\delta I)g,g\rangle
 \leq \langle(\Lambda-\Lambda_{\mathcal{B}}+\delta I)g,g\rangle + \delta \Vert g\Vert^2
 = (\theta+2\delta)\Vert g\Vert^2<0
\end{eqnarray*}
\noindent
for all $0<\delta<\delta_0:=-\frac{\theta}{2}$.
\\
On the other hand, if $\mathcal{B}\subseteq \mathcal{D}$, then for all $g\in L^2(\Gamma_{\mathrm{N}})^d$ 
\begin{eqnarray*}
\langle (\Lambda^\delta - \Lambda_{\mathcal{B}}+\delta I)g,g\rangle \geq -\delta \Vert g\Vert^2 + \delta \Vert g\Vert^2= 0. 
\end{eqnarray*}
\eproof
}
\\
\noindent
{\color{black}
We want to remark that analogous results to Lemma \ref{mono_test_D} and Corollary \ref{standard_test_noise} also hold for the case $\lambda\leq \lambda_0$ and $\mu\leq \mu_0$.
\\
\begin{remark}
For ill-posed inverse problems the term convergence rate usually denotes how fast the reconstruction converges to the 
true unknown as a function of the noise going to zero. Corollary \ref{standard_test_noise} shows that, for a fixed test inclusion and contrast
level, our criterion gives the correct answer when the noise is below a certain finite threshold. Hence, in this sense,
the theoretic convergence rate is infinitely fast. But in the numerical results we will see that the reconstructions
are in fact highly affected by noise. The reason is twofold. First, the required threshold for exactness may be low.
Second, in the discrete case, we test (\ref{test_standard_noise}) only for finitely many boundary loads $g_i$. If some eigenvectors to
negative eigenvalues are missing from the set $\{g_1,g_2,\ldots,g_m\}$, the quality of the test is reduced so that some blocks may be
wrongly assigned. However, this problem can be addressed by choosing $g_i$ as a sufficiently large basis of $L^2(\Gamma_{\mathrm{N}})^d$. 
It would be highly desirable to increase the method's noise robustness, e.g. as in \cite{Harrach_Mach}.
\end{remark}
}
\subsection{Numerical realization}

In order to close this subsection, we take a look at the numerical realization of the monotonicity test (Corollary \ref{standard_test_noise}) implemented with COMSOL Multiphysics
with LiveLink for MATLAB. 

\subsection*{Implementation}

We are given discrete noisy difference measurements $\overline{\Lambda}^\delta$ which fulfill
\begin{eqnarray}
 \Vert\overline{\Lambda}^\delta - \overline{\Lambda}_\mathcal{D}\Vert<\delta,
\end{eqnarray}
where the notation
\begin{eqnarray}
 \overline{A}:=(\langle g_i,A g_j\rangle)_{i,j=1}^m
\end{eqnarray}
\noindent
represents the discretized operator corresponding to $A\in\mathcal{L}(L^2(\Gamma_{\mathrm{N}})^d)$ w.r.t$.$ the boundary loads $g_i\in G$.
Here, the set $G:=\lbrace g_1, g_2,\ldots, g_m\rbrace$ is a system of boundary loads $g_i$, in which the boundary loads $g_i$ and $g_j$ ($i\neq j$) are pairwise orthogonal.
\\
\\
Let $(\lambda_0,\mu_0)$
be the Lam\'e parameters of the background material and $(\lambda_1,\mu_1)$ be the Lam\'e parameters of the inclusion,
so that 
 \begin{eqnarray}\label{param_inclu_1}
  \lambda=
   \left \{ \begin{array}{ll}
  \lambda_0,&x\in\Omega\setminus\mathcal{D},\\
   \lambda_1,& x\in\mathcal{D},
  \end{array}\right.
 \\
  \mu=
 \left \{ \begin{array}{ll}
  \mu_0,& x\in \Omega\setminus\mathcal{D},\\
  \mu_1,& x\in\mathcal{D},
  \end{array}\right.\label{param_inclu_2}\\ 
 \tilde{\lambda}_k=
 \left \{ \begin{array}{ll}
  \lambda_0,& x\in \Omega\setminus \mathcal{B}_k,\\
  \lambda_0+\alpha, & x\in \mathcal{B}_k,
 \end{array}\right.\label{param_inclu_3}\\ 
 \tilde{\mu}_k=
 \left \{ \begin{array}{ll}
  \mu_0,& x\in \Omega\setminus \mathcal{B}_k,\\
  \mu_0+\beta,& x\in \mathcal{B}_k,
 \end{array}\right.\label{param_inclu_4}
 \end{eqnarray}
\noindent
where $\mathcal{D}$ denotes the inclusion to be detected and $\mathcal{B}_k$ are $k=1,...,n$ known test inclusions.
In addition, the contrasts must fulfill $\alpha\leq\lambda_1-\lambda_0$ and $\beta\leq \mu_1-\mu_0$. 
\\
\\
\noindent
For our numerical experiments, we simulate these discrete measurements by solving

\begin{equation}\hspace{-1cm}\left\{ \begin{array}{rcll}
\nabla\cdot\left(\lambda_0 (\nabla\cdot u_0)I + 2\mu_0 \hat{\nabla} u_0 \right)&=&0&\mathrm{in}\,\,\Omega,\label{monotonicity_test_1}\\
-\nabla\cdot\left(((\lambda_1-\lambda_0)\chi_D) (\nabla\cdot u_0)I + 2((\mu_1-\mu_0)\chi_D) \hat{\nabla} u_0 \right)\\ \nonumber
+\nabla\cdot \left(\lambda (\nabla\cdot v)I + 2\mu \hat{\nabla} v \right)&=&0 &\mathrm{in}\,\,\Omega,\\
\left(\lambda_0 (\nabla\cdot u_0)I + 2\mu_0 \hat{\nabla} u_0 \right) n&=&g_i&\mathrm{on}\,\, \Gamma_{\mathrm{N}},\\
\left(\lambda (\nabla\cdot v)I + 2\mu \hat{\nabla} v \right) n&=&0&\mathrm{on}\,\, \Gamma_{\mathrm{N}},\\
u_0&=&0 &\mathrm{on}\,\, \Gamma_{\mathrm{D}},\\
v&=&0 &\mathrm{on}\,\, \Gamma_{\mathrm{D}},
\end{array}\right.
\end{equation}
\noindent
for each of the $i=1,\ldots,m$ boundary loads $g_i$, where $v:=u_0-u$. The equations regarding $v$ in the system (\ref{monotonicity_test_1}) 
result from substracting the boundary value problem (\ref{direct_1}) for the respective Lam\'{e} parameters.
\\
\\
To reconstruct the unknown inclusion $\mathcal{D}$, we determine the Neumann-to-Dirichlet operator 
for small test cubes $\mathcal{B}_k$, $k=1,\ldots n$, so that the
Lam\'e parameters are defined as in (\ref{param_inclu_3}) and (\ref{param_inclu_4}).
We calculate 
\begin{eqnarray*}
\overline{\Lambda}_k:=\overline{\Lambda}(\lambda_0,\mu_0)-\overline{\Lambda}(\tilde{\lambda}_k,\tilde{\mu}_k)
\end{eqnarray*}
\noindent
by solving an analogous system to (\ref{monotonicity_test_1}) (but with a different FEM
grid to avoid the so-called "inverse crime"). Note that this calculation does not depend on the measurements $\overline{\Lambda}^\delta$ and can be done in 
advance (in a so-called offline phase).
\noindent
We then compute the eigenvalues of
\begin{eqnarray*}
\overline{\Lambda}^\delta-\overline{\Lambda}_k+\delta I. 
\end{eqnarray*}
\noindent
If all eigenvalues are non-negative, then the test cube $\mathcal{B}_k$ is marked as "inside the inclusion".

\subsection*{Results}

\begin{figure}[H]
 \begin{center}
 \includegraphics[width=0.49\textwidth]{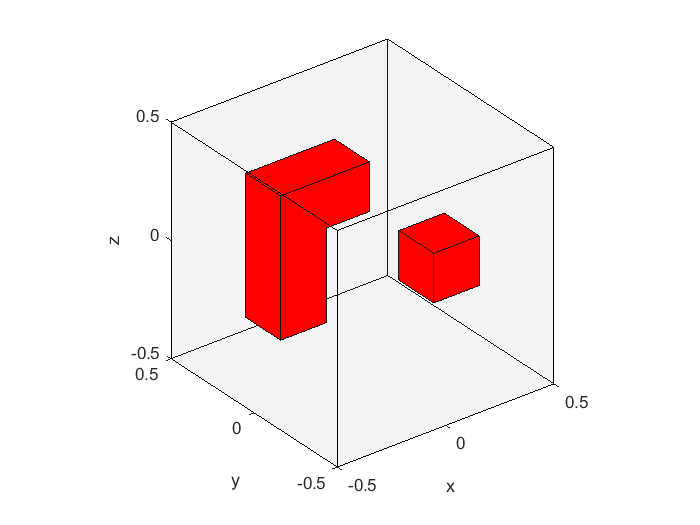}
  \includegraphics[width=0.49\textwidth]{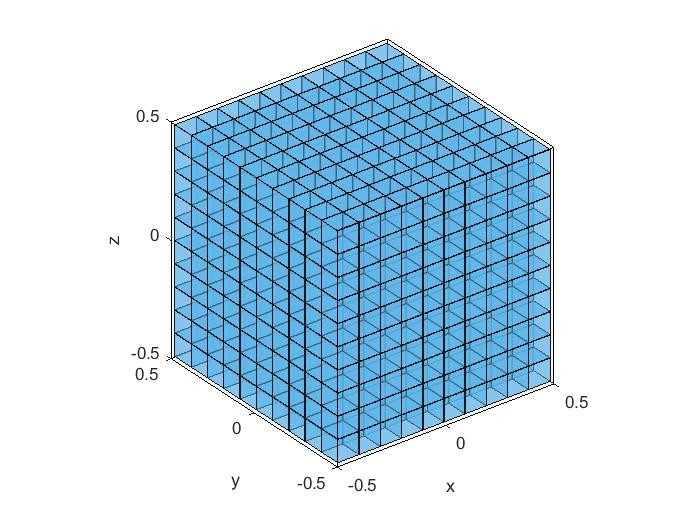}
\caption{Cube with two inclusions (red, left) and cube with $1000$ test blocks (blue, right).}\label{standard}
\end{center}
 \end{figure}

We present a test model, where we consider a cube of a biological tissue with two inclusions (tumors) as depicted in Figure \ref{standard}.
The Lam\'e parameters of the corresponding materials are given in Table \ref{lame_parameter_mono} (see \cite{Hiltawsky}).

 \begin{table} [H]
 \begin{center}

 \begin{tabular}{ |c|c| c |}  
\hline
 material & $\lambda$ & $\mu$ \\
  \hline
$x\in\Omega\setminus\mathcal{D}$: tissue &  $6.6211\cdot 10^5$   &  $6.6892\cdot 10^3$   \\
 \hline
\hspace{-0.5cm}$x\in\mathcal{D}$: tumor &  $2.3177\cdot 10^6$ &  $2.3411\cdot 10^4$  \\
\hline
\end{tabular}

\end{center}
\caption{Lam\'e parameter of the test material in [Pa].}
\label{lame_parameter_mono}
\end{table}
 \noindent
\noindent 
We use the $n=10\times 10 \times 10$ test cubes as shown in Figure \ref{standard}. The face characterized by $z=-0.5$
has zero displacement and each of the other faces is divided in $25$ squares of the same size resulting in $m=125$ patches, where the boundary loads $g_i$ are applied.
In this example, the boundary loads $g_i$ are the normal vectors on each patch. 

\begin{figure}[H]
 \begin{center}
 \includegraphics[width=0.49\textwidth]{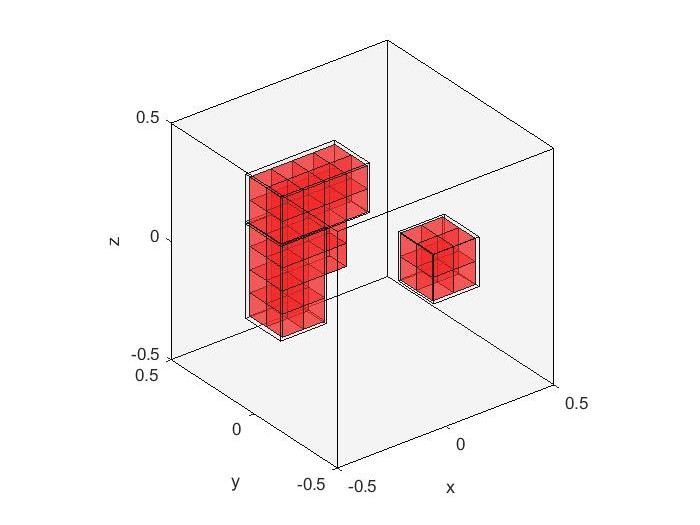}
\caption{Shape reconstruction of two inclusions (red) via $1000$ test cubes for $\alpha=(\lambda_1-\lambda_0)\approx 16.5\cdot 10^5$Pa, $\beta=(\mu_1-\mu_0) \approx 16.7\cdot 10^3$Pa without noise and $\delta=0$.}\label{result_1}
\end{center}
 \end{figure}

 \noindent
By performing the procedure described above without noise and using the parameters as given in Table \ref{appendix}, we end up with the result in Figure \ref{result_1}, 
depicting the test blocks with positive eigenvalues marked in red.

\begin{figure}[H]
 \begin{center}
 \includegraphics[width=0.49\textwidth]{geometry_new}
  \includegraphics[width=0.49\textwidth]{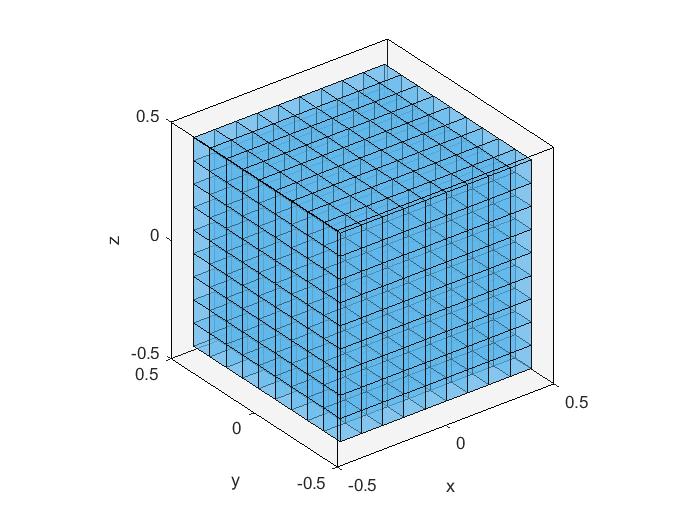}
\caption{Cube with two inclusions (red, left) and cube with $729$ test blocks (blue, right).}\label{standard_shift}
\end{center}
 \end{figure}

\noindent 
We can see that the inclusions are detected correctly and a clear separation of the two inclusions is obtained.
However, additional blocks were wrongly detected as lying inside the inclusions. 
\\
\\
Further on, we slightly change the setting and take a look at test cubes with the same size as in the first example (see Figure \ref{standard}).
However, now we shift these test cubes so that they do not ``perfectly fit'' into the unknown inclusions  as shown in Figure \ref{standard_shift}.

\begin{figure}[H]
 \begin{center}
 \includegraphics[width=0.49\textwidth]{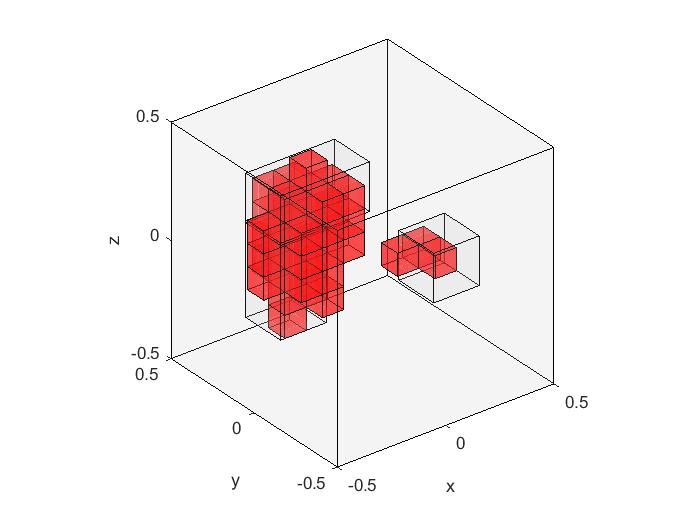}
  \includegraphics[width=0.49\textwidth]{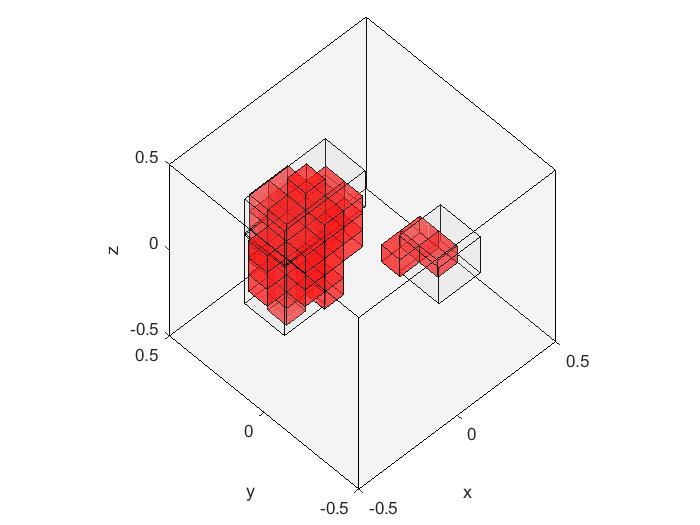}\\
  \includegraphics[width=0.55\textwidth]{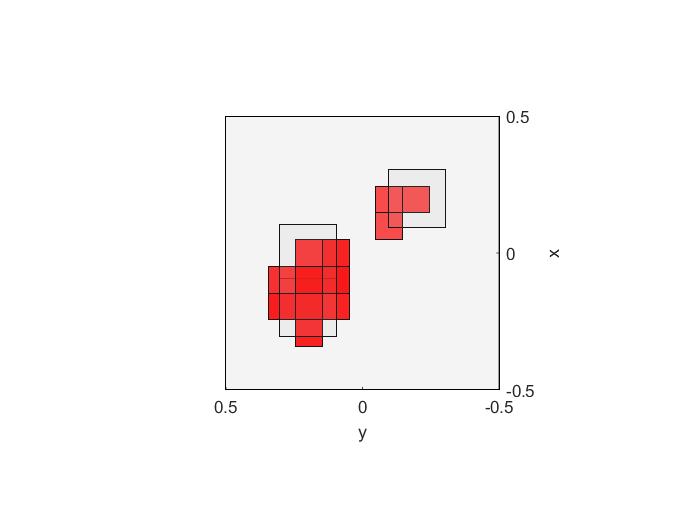} 
\caption{Shape reconstruction of two inclusions (red) via $729$ test cubes for $\alpha=(\lambda_1-\lambda_0)\approx 16.5\cdot 10^5$Pa, 
$\beta=(\mu_1-\mu_0) \approx 16.7\cdot 10^3$Pa without noise and $\delta=0$.}\label{result_stand_verschoben}
\end{center}
 \end{figure}

\begin{figure}[H]
 \begin{center}
 \includegraphics[width=0.49\textwidth]{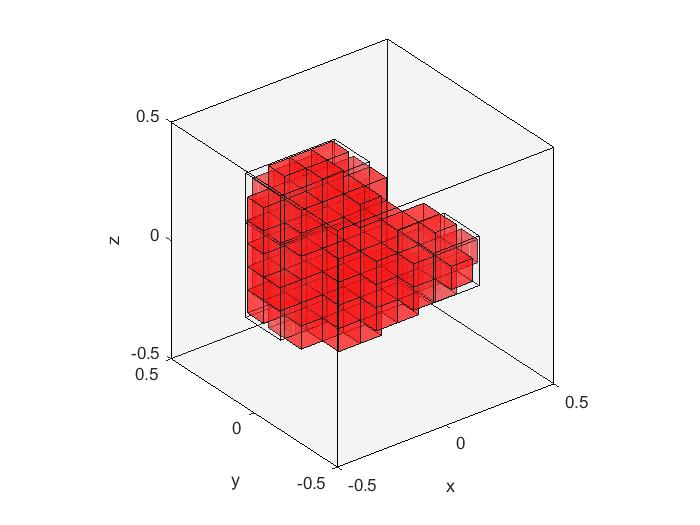}
\caption{Shape reconstruction of two inclusions (red) via $729$ test cubes
for $\alpha=(\lambda_1-\lambda_0)\approx 16.5\cdot 10^5$Pa, $\beta=(\mu_1-\mu_0) \approx 16.7\cdot 10^3$Pa
with a noise level of $0.1 \%$ and $\delta=3\cdot 10^{-10}$.}\label{result_2}
\end{center}
 \end{figure}

 \noindent
 Figure \ref{result_stand_verschoben} depicts the reconstruction of the two inclusions, which are detected but there are also blocks in between 
 wrongly marked.
\noindent
Next, we take a look at noisy measurement data $\overline{\Lambda}^{\delta}$ with a noise level of $0.1\%$ and set
\begin{eqnarray*}
\overline{\Lambda}^\delta= \overline{\Lambda}+10^{-3}(\overline{\Lambda}_{ij}E_{ij})_{i,j=1}^m, 
\end{eqnarray*}
\noindent
where $E\in \mathbb{R}^{m\times m}$ contains randomly distributed entries between $-1$ and $1$. 
{\color{black}
In doing so,
we obtain the reconstruction shown in Figure \ref{result_2}, where we chose $\delta=3\cdot 10^{-10}$ heuristically. 
In addition, we want to remark that this $\delta $ seems to be relatively small, if one observes that the absolute 
value of the largest eigenvalue is approximately of order $7\cdot 10^{-7}$. 
}

 \noindent
Thus, with a noise level of $0.1\%$ the two inclusions are still located correctly but they are not separated as two single inclusions due to noise.

\begin{remark}
The main drawback of the method is its offline computation time (see Table \ref{appendix} in the Appendix) due to the required solutions of problem (\ref{monotonicity_test_1})
for each test cube and each Neumann boundary load $g_i$. However, it should be noted that for testing multiple objects with the same background domain, we only have to 
perform
the offline phase of the test once, since it will be the same for each object to be tested. Hence, in this case, the computation time of the offline phase can be neglected. For
testing multiple objects with different background domains, we need a method whose offline computation time is faster. That method is the linearized monotonicity test discussed
in the next subsection.
\end{remark}

\subsection{Linearized monotonicity method}
Now, we go over to the modification of the standard monotonicity method considered so far and introduce the linearized version.
For that we apply the Fr\'echet derivative as defined in (\ref{Frechet}) in the form
\begin{eqnarray}\label{Frechet_lin}
\langle &\Lambda'(\lambda_0,\mu_0)(\alpha\chi_\mathcal{B},\beta \chi_\mathcal{B})g,h\rangle\\ \nonumber
&\quad=-\int_{\Omega} 2 \beta \chi_\mathcal{B}\, \hat{\nabla}u^{g}_{(\lambda_0,\mu_0)} : \hat{\nabla}u^{h}_{(\lambda_0,\mu_0)} 
+ \alpha\chi_\mathcal{B} \nabla \cdot u^{g}_{(\lambda_0,\mu_0)} \,\nabla\cdot  u^{h}_{(\lambda_0,\mu_0)}\,dx\\ \nonumber
&\quad=-\int_{\mathcal{B}} 2 \beta \, \hat{\nabla}u^{g}_{(\lambda_0,\mu_0)} : \hat{\nabla}u^{h}_{(\lambda_0,\mu_0)} 
+ \alpha \nabla \cdot u^{g}_{(\lambda_0,\mu_0)} \,\nabla\cdot  u^{h}_{(\lambda_0,\mu_0)}\,dx.
\end{eqnarray}

\begin{theorem}\label{theor_inclusion_2}
{\color{black}
Let $\lambda_0,\mu_0>0$ be constant and let $(\lambda,\mu)\in {\color{black}L^\infty_+(\Omega)\times L^\infty_+(\Omega)}$, 
$\lambda\geq \lambda_0$, $\mu\geq \mu_0$.}
For every open set $\mathcal{B}$ (e.g. ball or cube) and
every $\alpha,\beta\geq 0$, $\alpha+\beta>0$,
\begin{eqnarray*}
\alpha \chi_\mathcal{B}\leq\frac{\lambda_0}{\lambda} (\lambda -\lambda_0), \quad \beta \chi_\mathcal{B}\leq \frac{\mu_0}{\mu}(\mu -\mu_0)
\end{eqnarray*}
\noindent
implies
\begin{eqnarray*}
\Lambda(\lambda_0,\mu_0) +\Lambda^\prime(\lambda_0,\mu_0)(\alpha\chi_\mathcal{B},\beta\chi_\mathcal{B}) \geq \Lambda(\lambda,\mu)
\end{eqnarray*}
\noindent
and
\begin{eqnarray*}
\mathcal{B}\nsubseteq \mathcal{D}
\,\,\textnormal{\it implies}\,\,
\Lambda(\lambda_0,\mu_0)+\Lambda'(\lambda_0,\mu_0)(\alpha\chi_\mathcal{B},\beta\chi_\mathcal{B})\ngeq \Lambda(\lambda,\mu).
\end{eqnarray*}
Hence, the set
\begin{eqnarray*}
R:=\bigcup\limits_{\alpha,\beta\geq 0, \alpha+\beta>0}\lbrace \mathcal{B} \subseteq\Omega: \Lambda(\lambda_0,\mu_0)+\Lambda^\prime(\lambda_0,\mu_0) (\alpha\chi_\mathcal{B},\beta\chi_\mathcal{B})\geq \Lambda(\lambda,\mu)\rbrace
\end{eqnarray*}
\noindent
fulfills
\begin{eqnarray*}
\mathrm{inn}\,\mathrm{supp}((\lambda -\lambda_0, \mu -\mu_0)^T) \subseteq R \subseteq 
\mathcal{D}.
\end{eqnarray*}
\end{theorem}

\bproof
{\color{black}
Let $\lambda_0,\mu_0>0$ be constant and let $(\lambda,\mu)\in {\color{black}L^\infty_+(\Omega)\times L^\infty_+(\Omega)}$, 
$\lambda\geq \lambda_0$, $\mu\geq \mu_0$.} Let $\mathcal{B}$ be an open set and 
$\alpha,\beta\geq 0$, $\alpha+\beta>0$.
\\
{\color{black}
For every $g\in L^2(\Gamma_{\mathrm{N}})^d$ and solution $u_{(\lambda_0,\mu_0)}^g\in \mathcal{V}$ of
\begin{equation}\left\{ \begin{array}{rcll}\label{theo_4}
\nabla\cdot \left(\lambda_0 (\nabla\cdot u_{(\lambda_0,\mu_0)})I + 2 \mu_0 \hat{\nabla} u_{(\lambda_0,\mu_0)}^g \right)&=&0 &\quad \mathrm{in}\,\,\Omega,\\
\left(\lambda_0 (\nabla\cdot u_{(\lambda_0,\mu_0)}^g)I + 2 \mu_0 \hat{\nabla} u_{(\lambda_0,\mu_0)}^g \right) n&=&g&\quad \mathrm{on}\,\, \Gamma_{\mathrm{N}},\\
u_{(\lambda_0,\mu_0)}^g&=&0 &\quad \mathrm{on}\,\, \Gamma_{\mathrm{D}},
\end{array}\right.
\end{equation}
}
\noindent
we have with Lemma \ref{mon_est_2} 
\begin{eqnarray*}
 &\langle g, (\Lambda(\lambda_0,\mu_0)+\Lambda^\prime(\lambda_0,\mu_0)(\alpha\chi_\mathcal{B},\beta\chi_\mathcal{B})-\Lambda(\lambda,\mu))g\rangle\\
 &\quad\geq  \int_\Omega\Bigg( 2\left(\frac{\mu_0}{\mu}(\mu-\mu_0)-\beta\chi_\mathcal{B}\right) \hat{\nabla}u_{(\lambda_0,\mu_0)}^g:\hat{\nabla}u_{(\lambda_0,\mu_0)}^g\\
 &\qquad+\left(\frac{\lambda_0}{\lambda}(\lambda-\lambda_0)-\alpha\chi_\mathcal{B}\right)\nabla\cdot u_{(\lambda_0,\mu_0)}^g\nabla\cdot u_{(\lambda_0,\mu_0)}^g\Bigg)\,dx.
\end{eqnarray*}
\noindent
\\
Hence, $\alpha\chi_\mathcal{B}\leq \frac{\lambda_0}{\lambda}(\lambda-\lambda_0)$, $\beta\chi_\mathcal{B}\leq \frac{\mu_0}{\mu}(\mu-\mu_0)$ implies
\begin{eqnarray*}
 \Lambda(\lambda_0,\mu_0)+\Lambda^\prime(\lambda_0,\mu_0)(\alpha\chi_\mathcal{B},\beta\chi_\mathcal{B})\geq \Lambda(\lambda,\mu).
\end{eqnarray*}
\noindent
It remains to show that
\begin{eqnarray}\label{lin_mon_B_not_in_D}
\mathcal{B}\nsubseteq \mathcal{D}\quad \mathrm{implies}\quad
\Lambda(\lambda_0,\mu_0)+\Lambda^\prime(\lambda_0,\mu_0)(\alpha\chi_\mathcal{B}, \beta\chi_\mathcal{B} )\ngeq \Lambda(\lambda,\mu).
\end{eqnarray}
{\color{black}
Let $\mathcal{B}\nsubseteq \mathcal{D}$.
The  monotonicity relation from Corollary \ref{monotonicity} yields that shrinking the open set $\mathcal{B}$ makes
$\Lambda(\lambda_0,\mu_0)+\Lambda^\prime(\lambda_0,\mu_0)(\alpha\chi_\mathcal{B}, \beta\chi_\mathcal{B} )$ larger, so  we can assume without loss of generality that 
$\mathcal{B}\subset \Omega \setminus \mathcal{D}$.
Then
\begin{eqnarray*}
&\fl\langle g,(\Lambda(\lambda_0,\mu_0)+\Lambda'(\lambda_0,\mu_0)(\alpha\chi_{\mathcal{B}}, \beta\chi_{\mathcal{B}} )-\Lambda(\lambda,\mu)g\rangle\\
&\fl\quad\leq \int_\Omega \bigg( 2 (\mu-\mu_0-\beta\chi_{\mathcal{B}}) \hat{\nabla}u_{(\lambda_0,\mu_0)}^g : \hat{\nabla}u_{(\lambda_0,\mu_0)}^g+
(\lambda-\lambda_0-\alpha\chi_{\mathcal{B}})\nabla\cdot u_{(\lambda_0,\mu_0)}^g \nabla\cdot u_{(\lambda_0,\mu_0)}^g\bigg)\,dx\\
&\fl\quad\leq \int_{\mathcal{D}} \bigg( 2 (\mu-\mu_0) \hat{\nabla}u_{(\lambda_0,\mu_0)}^g : \hat{\nabla}u_{(\lambda_0,\mu_0)}^g+
(\lambda-\lambda_0)\nabla\cdot u_{(\lambda_0,\mu_0)}^g \nabla\cdot u_{(\lambda_0,\mu_0)}^g\bigg)\,dx\\
 &\fl\qquad-\int_{\mathcal{B}}\left(2\beta\hat{\nabla}u_{(\lambda_0,\mu_0)}^g : \hat{\nabla}u_{(\lambda_0,\mu_0)}^g+
 \alpha \nabla\cdot u_{(\lambda_0,\mu_0)}^g \nabla\cdot u_{(\lambda_0,\mu_0)}^g\right)\,dx\\
 \end{eqnarray*}
\noindent
so that the assertion follows using localized potential for the background  parameters  $(\mu_0,\lambda_0)$  as in Theorem \ref{monotonicity}.
}
\eproof
\noindent
\\
\\
Next, we formulate a theorem for the case $\lambda\leq \lambda_0$, $\mu\leq \mu_0$.
{\color{black}
\begin{theorem}\label{theor_inclusion_lin_add}
{\color{black}
Let $\lambda_0,\mu_0>0$ be constant and let $(\lambda,\mu)\in {\color{black}L^\infty_+(\Omega)\times L^\infty_+(\Omega)}$, 
$\lambda\leq \lambda_0$, $\mu\leq \mu_0$.}
For every open set $\mathcal{B}$ (e.g. ball or cube) and
every $\alpha,\beta\geq 0$, $\alpha+\beta>0$,
\begin{eqnarray*}
\alpha \chi_\mathcal{B}\leq \lambda_0 -\lambda, \quad \beta \chi_\mathcal{B}\leq\mu_0 -\mu
\end{eqnarray*}
\noindent
implies
\begin{eqnarray*}
\Lambda(\lambda_0,\mu_0) -\Lambda^\prime(\lambda_0,\mu_0)(\alpha\chi_\mathcal{B},\beta\chi_\mathcal{B}) \leq \Lambda(\lambda,\mu)
\end{eqnarray*}
\noindent
and
\begin{eqnarray*}
\mathcal{B}\nsubseteq \mathcal{D}
\,\,\textnormal{\it implies}\,\,
\Lambda(\lambda_0,\mu_0)-\Lambda'(\lambda_0,\mu_0)(\alpha\chi_\mathcal{B},\beta\chi_\mathcal{B})\nleq \Lambda(\lambda,\mu).
\end{eqnarray*}
Hence, the set
\begin{eqnarray*}
R:=\bigcup\limits_{\alpha,\beta\geq 0, \alpha+\beta>0}\lbrace \mathcal{B} \subseteq\Omega: \Lambda(\lambda_0,\mu_0)-\Lambda^\prime(\lambda_0,\mu_0) (\alpha\chi_\mathcal{B},\beta\chi_\mathcal{B})\leq \Lambda(\lambda,\mu)\rbrace
\end{eqnarray*}
\noindent
fulfills
\begin{eqnarray*}
\mathrm{inn}\,\mathrm{supp}((\lambda_0-\lambda,\mu_0-\mu)^T)\subseteq R \subseteq 
\mathcal{D}.
\end{eqnarray*}
\end{theorem}

\bproof
In order to prove this theorem, we have to consider the same steps as in the proof of Theorem \ref{theor_inclusion_2}, 
but we have to apply the estimate from Lemma \ref{mono} instead of Lemma \ref{mon_est_2} which results in different upper bounds on the contrasts $\alpha$ and $\beta$.
\eproof
\noindent
\\
\\
Summarizing the results from Theorem \ref{theor_inclusion_2} as well as Theorem \ref{theor_inclusion_lin_add}, we end up with the linearized monotonicity tests as introduced in Corollary 
\ref{cor_lin_mon_1} and \ref{cor_lin_mon_2}. 
}

\begin{corollary}{Linearized monotonicity test: 1. version}\label{cor_lin_mon_1}
\\
{\color{black}
Let $\lambda_0$, $\lambda_1$, $\mu_0$, $\mu_1\in\mathbb{R}^+$ with $\lambda_1>\lambda_0$, $\mu_1>\mu_0$  
and assume that 
\mbox{$(\lambda,\mu)=$} \mbox{$(\lambda_0+(\lambda_1-\lambda_0)\chi_\mathcal{D},\mu_0+(\mu_1-\mu_0)\chi_{\mathcal{D}})$}
with $\mathcal{D}=\mathrm{out}_{\partial\Omega}\,\mathrm{supp}((\lambda-\lambda_0,\mu-\mu_0)^T)$.}
Further on, let $\alpha,\beta\geq 0$, $\alpha+\beta>0$ with $\alpha \leq \frac{\lambda_0}{\lambda_1} (\lambda_1 -\lambda_0) $,
$\beta\leq \frac{\mu_0}{\mu_1} (\mu_1 -\mu_0) $.
Then for every open set $\mathcal{B}\subseteq \Omega$
\begin{eqnarray*}
\mathcal{B}\subseteq\mathcal{D}\quad\mathrm{if\,and\,only\,if}\quad \Lambda(\lambda_0,\mu_0)+\Lambda^\prime(\lambda_0,\mu_0)\left(\alpha\chi_\mathcal{B}, \beta\chi_\mathcal{B}\right)\geq \Lambda(\lambda,\mu).
\end{eqnarray*}
\end{corollary}
{\color{black}
\bproof
The proof of Corollary \ref{cor_lin_mon_1} follows from Theorem \ref{theor_inclusion_2} analogous to the proof of Corollary \ref{col_stan_mon_1}.
\eproof
}
\noindent
\\
\\
In addition, we formulate the linearized monotonicity test for the case $\lambda\leq \lambda_0$ and $\mu\leq \mu_0$.
\begin{corollary}{Linearized monotonicity test: 2. version}\label{cor_lin_mon_2}
\\
{\color{black}
Let $\lambda_0$, $\lambda_1$, $\mu_0$, $\mu_1\in\mathbb{R}^+$ with $\lambda_1<\lambda_0$, $\mu_1<\mu_0$  
and assume that 
\mbox{$(\lambda,\mu)=$} \mbox{$(\lambda_0+(\lambda_1-\lambda_0)\chi_\mathcal{D},\mu_0+(\mu_1-\mu_0)\chi_{\mathcal{D}})$}
with $\mathcal{D}=\mathrm{out}_{\partial\Omega}\,\mathrm{supp}((\lambda-\lambda_0,\mu-\mu_0)^T)$.}
Further on, let $\alpha,\beta\geq 0$, $\alpha+\beta>0$ with $\alpha\leq\lambda_0 -\lambda_1 $,
$ \beta\leq \mu_0 -\mu_1 $.
Then for every open set $\mathcal{B}\subseteq \Omega$
\begin{eqnarray*}
\mathcal{B}\subseteq\mathcal{D}\quad\mathrm{if\,and\,only\,if}\quad \Lambda(\lambda_0,\mu_0)-\Lambda^\prime(\lambda_0,\mu_0)\left(\alpha\chi_\mathcal{B},\beta \chi_\mathcal{B}\right)\leq \Lambda(\lambda,\mu).
\end{eqnarray*}
\end{corollary}
{\color{black}
\bproof
In order to prove this corollary, we apply the results of Theorem \ref{theor_inclusion_lin_add} in a similar way as in the the proof of Corollary \ref{col_stan_mon_1}.
\eproof
}
\noindent
\\
\\
Next, we apply Theorem \ref{theor_inclusion_2} to difference measurements $\Lambda_\mathcal{D}$ as defined in (\ref{Lambda_D}) and obtain the following lemma.

\begin{lemma}\label{lin_mono_test_delta}
Let $\Lambda^\prime_\mathcal{B}=\Lambda^\prime(\lambda_0,\mu_0)(\alpha\chi_\mathcal{B},\beta\chi_\mathcal{B})$. Under the same assumptions on $\lambda$ and $\mu$
as in Theorem \ref{theor_inclusion_2}, we have for every open set $\mathcal{B}$
(e.g. ball or cube) and every $\alpha,\beta\geq 0$, $\alpha+\beta>0$,
\noindent
\begin{eqnarray}\label{alpha_beta_lin}
\alpha \chi_\mathcal{B}\leq\frac{\lambda_0}{\lambda} (\lambda -\lambda_0), \quad \beta \chi_\mathcal{B}\leq \frac{\mu_0}{\mu}(\mu -\mu_0)
\end{eqnarray}
implies
\begin{eqnarray*}
\Lambda_\mathcal{D}+\Lambda^\prime_\mathcal{B} \geq 0
\end{eqnarray*}
\noindent
and
\begin{eqnarray*}
\mathcal{B}\nsubseteq\mathcal{D}
\,\,\textnormal{\it implies}\,\,
\Lambda_\mathcal{D} +\Lambda^\prime_\mathcal{B} \ngeq 0.
\end{eqnarray*}
\end{lemma}
{\color{black}
\bproof
Similar as in the proof of Theorem \ref{theor_inclusion_1}, we start with the consideration of
$\Lambda_\mathcal{D}+\Lambda_\mathcal{B}^\prime=\Lambda(\lambda_0,\mu_0)-\Lambda(\lambda,\mu)+\Lambda^\prime(\lambda_0,\mu_0)(\alpha\chi_\mathcal{B},\beta\chi_\mathcal{B})$
and obtain the corresponding relations with Theorem \ref{theor_inclusion_2}.

\eproof
}
\noindent
\\
As for the standard monotonicity test, we also introduce a linearized monotonicity test for noisy data (c.f$.$ Corollary \ref{cor_lin_mon_1}). Let the 
noisy difference data be given by (\ref{Lambda_delta}) and fulfill (\ref{Lambda_delta_norm}).
{\color{black}
\begin{corollary}{Linearized monotonicity test for noisy data}\label{lin_mono_test_noise}
\\
{\color{black}
Let $\lambda_0$, $\lambda_1$, $\mu_0$, $\mu_1\in\mathbb{R}^+$ with $\lambda_1>\lambda_0$, $\mu_1>\mu_0$  
and assume that 
\mbox{$(\lambda,\mu)=$} \mbox{$(\lambda_0+(\lambda_1-\lambda_0)\chi_\mathcal{D},\mu_0+(\mu_1-\mu_0)\chi_{\mathcal{D}})$}
with $\mathcal{D}=\mathrm{out}_{\partial\Omega}\,\mathrm{supp}((\lambda-\lambda_0,\mu-\mu_0)^T)$.}
Further on, let $\alpha,\beta\geq 0$, $\alpha+\beta>0$ with $\alpha \leq \frac{\lambda_0}{\lambda_1} (\lambda_1 -\lambda_0) $,
$\beta\leq \frac{\mu_0}{\mu_1} (\mu_1 -\mu_0) $.
{\color{black}
Let $\Lambda^\delta$ be the Neumann-to-Dirichlet operator for noisy difference measurements with noise level $\delta>0$.
Then for every open set $\mathcal{B}\subseteq\Omega$ there exists a noise level $\delta_0>0$, such that for all
$0<\delta<\delta_0$, $\mathcal{B}$ is correctly detected as inside or not inside the inclusion $\mathcal{D}$ 
by the following monotonicity test
\begin{eqnarray*}
\mathcal{B}\subseteq\mathcal{D}\quad \textnormal{\it if and only if} \quad \Lambda^\delta +\Lambda^\prime_\mathcal{B}+\delta I \geq 0.
\end{eqnarray*}
}
\end{corollary}
\bproof
We proceed in a similar manner as in the proof of Corollary \ref{standard_test_noise} and summarize the essential results of that proof. 
Let $\mathcal{B}\not\subseteq \mathcal{D}$. By Corollary \ref{compact} $\Lambda_{\mathcal{D}}+\Lambda_{\mathcal{B}}^\prime$
is compact and self-adjoint and by Lemma \ref{lin_mono_test_noise} $\Lambda_{\mathcal{D}}+\Lambda_{\mathcal{B}}^\prime\not\geq 0$.
Hence, $\Lambda_{\mathcal{D}}-\Lambda_{\mathcal{B}}$ has negative eigenvalues. Let $\theta<0$ be the smallest eigenvalue with corresponding
eigenvector $g\in L^2(\Gamma_{\mathrm{N}})^d$, so that
\begin{eqnarray*}
 \langle (\Lambda^\delta+\Lambda_{\mathcal{B}}^\prime+\delta I)g,g\rangle
 \leq \langle(\Lambda+\Lambda_{\mathcal{B}}^\prime+\delta I)g,g\rangle + \delta \Vert g\Vert^2
 = (\theta+2\delta)\Vert g\Vert^2<0
\end{eqnarray*}
\noindent
for all $0<\delta<\delta_0:=-\frac{\theta}{2}$.
\\
On the other hand, if $\mathcal{B}\subseteq \mathcal{D}$, then for all $g\in L^2(\Gamma_{\mathrm{N}})^d$ 
\begin{eqnarray*}
\langle (\Lambda^\delta + \Lambda_{\mathcal{B}}^\prime+\delta I)g,g\rangle \geq -\delta \Vert g\Vert^2 + \delta \Vert g\Vert^2= 0. 
\end{eqnarray*}
\eproof
}
\\
\noindent
{\color{black}
We want to remark that analogous results to Lemma \ref{lin_mono_test_delta} and Corollary \ref{lin_mono_test_noise} also hold for the case $\lambda\leq \lambda_0$ and $\mu\leq \mu_0$.
}

\subsection{Numerical realization}
Next, we go over to the implementation of the linearized monotonicity test according to Corollary \ref{lin_mono_test_noise}.
The implementation is again conducted  with COMSOL Multiphysics
with LiveLink for MATLAB.

\subsection*{Implementation}
We solve system (\ref{monotonicity_test_1}) to simulate our discrete noisy difference measurements $\overline{\Lambda}^\delta$
which must fulfill
\begin{eqnarray*}
 \Vert\overline{\Lambda}^\delta - \overline{\Lambda}\Vert<\delta.
\end{eqnarray*}
\noindent
With $\lambda$, $\mu$ as in (\ref{param_inclu_1})-(\ref{param_inclu_2}) and
\begin{eqnarray*}
 \tilde{\lambda}_k&=\alpha \chi_{\mathcal{B}_k},\\
 \tilde{\mu}_k&=\beta \chi_{\mathcal{B}_k},
\end{eqnarray*}
\noindent
where the contrasts must fulfill $\alpha\leq \frac{\lambda_0}{\lambda_1}(\lambda_1-\lambda_0)$ and
 $\beta\leq \frac{\mu_0}{\mu_1}(\mu_1-\mu_0)$, we proceed as follows.
To reconstruct the unknown inclusion $\mathcal{D}$, we determine the Fr\'echet derivative
\begin{eqnarray*}
\overline{\Lambda}^\prime_k:=\overline{\Lambda}^\prime(\lambda_0,\mu_0)(\tilde{\lambda}_k,\tilde{\mu}_k)
\end{eqnarray*}
\noindent
for each of the test cubes $\mathcal{B}_k$, $k=1,\ldots,n$, by an approximation via Gaussian quadrature in MATLAB.
{\color{black} For more details concerning Gaussian quadrature, the reader is referred to, e.g., \cite{Jinyun}.}
The required solution $u_0$ for 
the background Lam\'{e} parameters used in the calculation of $\overline{\Lambda}^\prime_k$ is again calculated via COMSOL.
\\
\\
Note that this calculation does not depend on the measurements $\Lambda^\delta$ and can be done in 
advance (in a so-called offline phase).
\\
\\
We then estimate the eigenvalues of
\begin{eqnarray*}
\overline{\Lambda}^\delta+\overline{\Lambda}^\prime_k+\delta I. 
\end{eqnarray*}
\noindent
If all eigenvalues are non-negative, then the test cube $\mathcal{B}_k$ is marked as "inside the inclusion".

\subsection*{Results}

\begin{figure}[H]
 \begin{center}
 \includegraphics[width=0.49\textwidth]{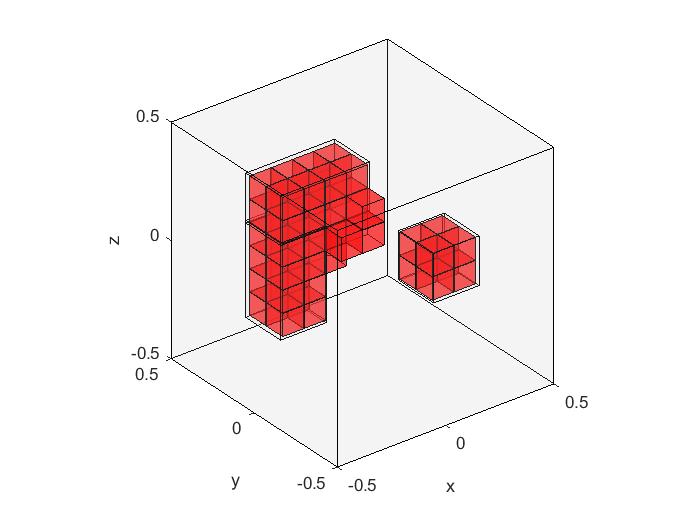}
\caption{Shape reconstruction of two inclusions (red) via $1000$ test cubes for
$\alpha=0.28(\lambda_1-\lambda_0)\approx 4.6\cdot 10^5$Pa, $\beta=0.28(\mu_1-\mu_0)\approx 4.7\cdot 10^3$ Pa 
without noise and $\delta=0$.}\label{result_lin_1000}
\end{center}
 \end{figure}

First, we consider the same setting as for the standard monotonicity method.
Implementing the linearized monotonicity method for the same boundary 
loads $g_i$ and without noise leads to the result depicted in Figure \ref{result_lin_1000},
where the test cubes with positive eigenvalues are marked in red.
The computation time as well as the applied parameters can again be found in the Appendix 
(Table \ref{appendix}). 

\noindent
As a result, we get a similar outcome as we obtained with the non-linearized monotonicity tests (for comparison see Figure \ref{result_1}).
However, the result for the non-linearized monotonicity method is slightly better. The reason is the stronger constraint on $\alpha$ and $\beta$
in the linearized case (c.f. Equations (\ref{alpha_beta}) and (\ref{alpha_beta_lin})). In other words, the non-linearized method allows us
to test with stronger contrasts.

\noindent

\begin{figure}[H]
 \begin{center}
 \includegraphics[width=0.49\textwidth]{geometry_new}
  \includegraphics[width=0.49\textwidth]{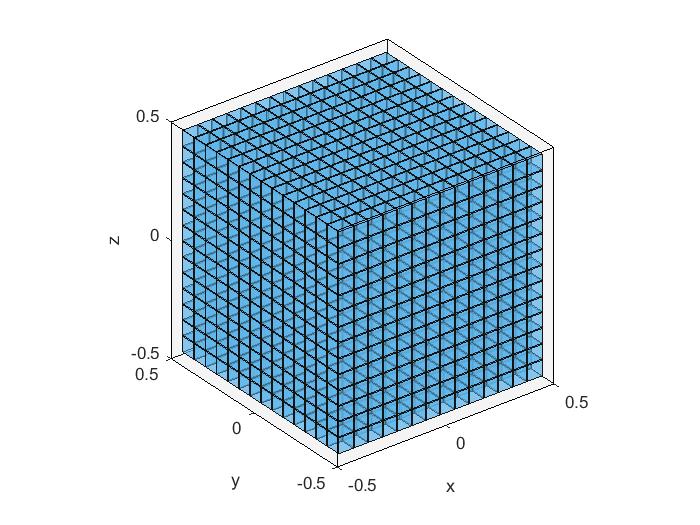}
\caption{Cube with two inclusions (red, left) and cube with $2744$ test blocks (blue, right).}\label{linearized}
\end{center}
 \end{figure}

\begin{remark}\label{adv_lin_mon}
All in all, the complete computation time including the offline phase for the linearized monotonicity test is faster (cf. Table \ref{appendix}).
However, the direct and the online time are similar if one uses the same mesh size. Hence, for testing objects with different background domains, the linearized
approach should be applied since the offline phase has to be realized each time. For testing many objects with the same background domain,
the non-linearized method should be applied, since it provides a slightly more precise reconstruction. The drawback of the offline calculation
time is mitigated in this case, since it has to be calculated only once at the beginning, so that the testing itself for each object is as fast as the
testing via the linearized method without its offline phase.
\end{remark}
\noindent
The significantly shorter offline phase of the linearized monotonicity method (see Table \ref{appendix}),
allows us to increase the number of test blocks. In addition, we also increase the number of tetrahedrons to test the influence of the
mesh size on the runtime and in order to provide enough nodes inside of the smaller test blocks so that the Fr\'echet derivative
can be calculated with sufficient precision. 
Hence, we go over to a slightly different setting (Figure \ref{linearized}), where we perform the testing
 for $125$ patches with $14\times 14\times 14$ test cubes instead 
 of $10\times10 \times 10$ cubes in order to obtain a higher resolution. 
We want to remark that these smaller cubes are chosen in such a way that some of them are also lying on the boundary of the inclusions we want to detect.
 That means, that in theory, we are able to detect also smaller inclusions than the test cubes used before
 and a higher resolution enables us to reconstruct inclusions, which do not surround a complete test block for a coarser set of test blocks, since
 Theorem \ref{theor_inclusion_2} guarantees that we cannot detect inclusions smaller than the used test inclusions. 
 
 \noindent
 
\begin{figure}[H]
 \begin{center}
 \includegraphics[width=0.49\textwidth]{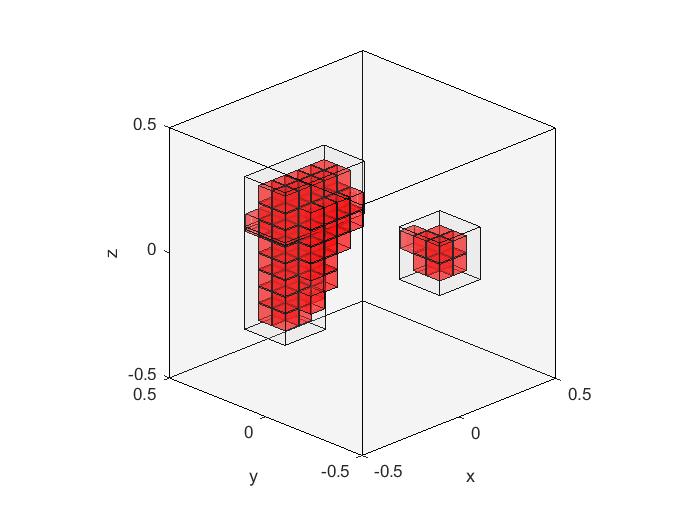}
  \includegraphics[width=0.49\textwidth]{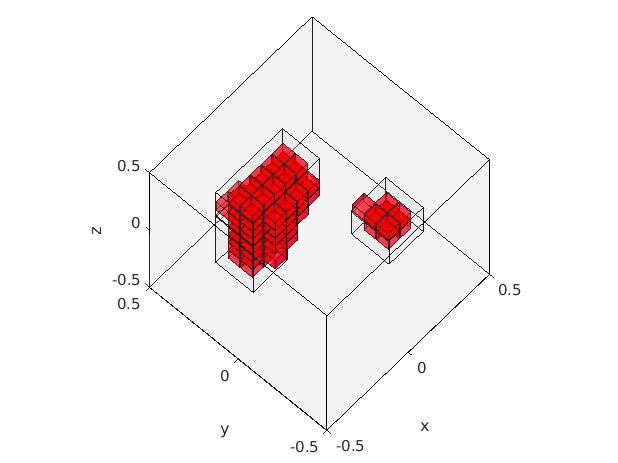}\\
   \includegraphics[width=0.55\textwidth]{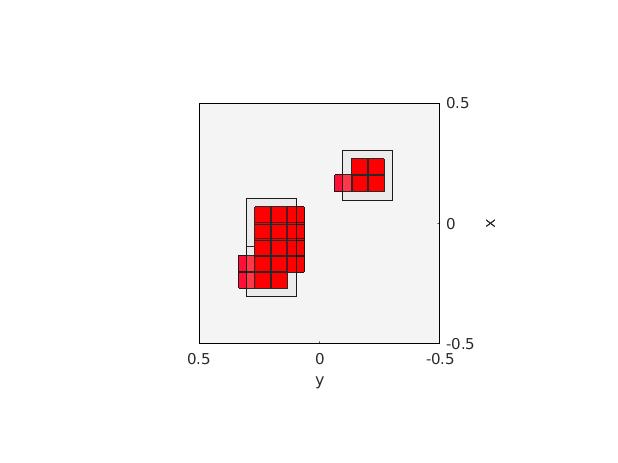}
\caption{Shape reconstruction of two inclusions (red) via $2744$ test cubes
for $\alpha=0.28(\lambda_1-\lambda_0)\approx 4.6\cdot 10^5$Pa, $\beta=0.28(\mu_1-\mu_0)\approx 4.7\cdot 10^3$ Pa
without noise and $\delta=0$.}\label{result_lin}
\end{center}
 \end{figure}

 \noindent
We see in Figure \ref{result_lin}, that the two inclusions are detected and their shape is reconstructed almost correctly.
However, additional blocks where wrongly detected. Those blocks lie again between the two inclusions,
as was the case for the standard monotonicity method, or are only partially located inside of the inclusion.
 \begin{remark}
  If we compare the noiseless results of the standard monotonicity method (Figure \ref{result_stand_verschoben}) with the one from the linearized monotonicity method 
 (Figure \ref{result_lin}), we can conclude that in our examples with both methods, the inclusions were detected. In addition,
 our numerical example showed that both monotonicity methods reconstruct the shape of the inclusions and are able to separate them.
 \end{remark}

\begin{figure}[H]
 \begin{center}
 \includegraphics[width=0.49\textwidth]{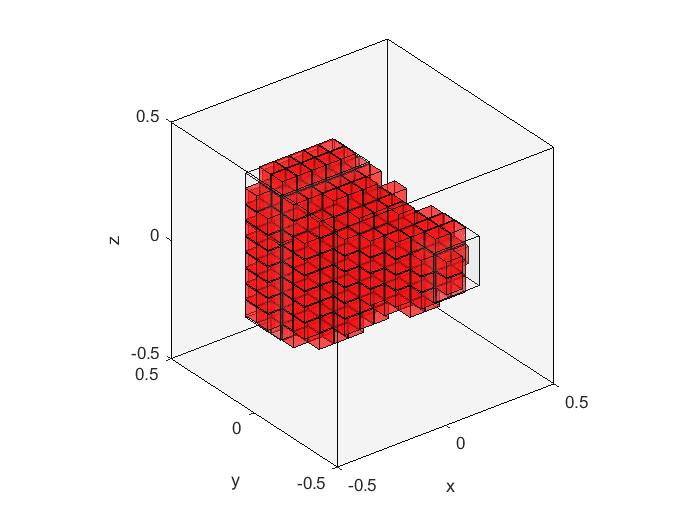}
\caption{Shape reconstruction of two inclusions (red) via $2744$ test cubes for
$\alpha=0.28(\lambda_1-\lambda_0)\approx 4.6\cdot 10^5$Pa, $\beta=0.28(\mu_1-\mu_0)\approx 4.7\cdot 10^3$ Pa with noise level of $0.1\%$ and $\delta=5\cdot10^{-10}$.}\label{result_lin_noise}
\end{center}
 \end{figure}

 \noindent
 Finally, as for the standard monotonicity method, we again present the result obtained from noisy data $\Lambda^{\delta}$ with a noise
 level of $0.1\%$ (Figure \ref{result_lin_noise}). 
 
\noindent
This reconstruction shows us that the linearized monotonicity method applied to noisy data gives similar results as the standard monotonicity method (see Figure \ref{result_2}).
The inclusion is again located correctly, but a separation is not possible due to the noise level.

\section{Conclusion and Outlook}

In this paper we introduced and analyzed the standard as well as linearized monotonicity method for the linear elasticity problem and gave an insight into
the performance of the monotonicity tests.
Further on, we showed numerical examples of the different methods and compared them with each other. The next step will be 
the adaptation of these methods to the elastic Helmholtz-type equation.

\section*{Appendix}

The runtime of the two methods are given in Table 2. We distinguish between the
solution of the direct and inverse problem, which itself is split into its offline and online phase.
The direct problem represents the simulation of the difference measurements, which are provided by 
a user in practice. For the inverse problem, the offline phase includes the calculation of all
quantities which can be computed without knowing the measurements, i.e. all quantities that 
only depend on $(\lambda_0,\mu_0)$ and are calculated on a different grid as the direct problem
to avoid ``inverse crime''. Contrary, the online phase includes all computations
which depend on the difference measurements.
\\
\\
It should be noted that the approach for the direct problem is the same for the non-linearized
and linearized method, so that the computation times are nearly identical for the same grid size (see the two examples for $1000$ pixel).
As expected, the computation time increases while increasing the number of tetrahedrons used in the discretization (cf. last line of Table \ref{appendix}).
In addition, the offline phase is considerably larger than the runtime of the direct problem.
This is due to the computation time of the Fr\'echet derivative.

\begin{table}[H]
\begin{center}
 \begin{tabular}{ |p{6.4cm}|p{2.6cm}|p{2.8cm}|p{2.8cm}|}  
\hline
runtime & direct problem & inverse problem $\quad$ $\quad$ offline phase & inverse problem $\quad$ $\quad$ online phase\\
\hline
non-linearized method, \newline$\alpha=(\lambda_1-\lambda_0)\approx 16.5\cdot 10^5$Pa,\newline $\beta=(\mu_1-\mu_0) \approx 16.7\cdot 10^3$Pa,\newline $1000$ pixel, $\#$tetrahedrons $=7644$ (direct problem, offline similar)& $9$min $10$s & $6$d $23$h $43$min & $2$min $27$s\\
\hline
non-linearized method, \newline$\alpha=(\lambda_1-\lambda_0)\approx 16.5\cdot 10^5$Pa,\newline $\beta=(\mu_1-\mu_0) \approx 16.7\cdot 10^3$Pa,\newline $729$ pixel, $\#$tetrahedrons $=7644$ (direct problem, offline similar)& $9$min $23$s & $4$d $20$h $54$min & $2$min $37$s\\
\hline
linearized method, \newline$\alpha=0.28(\lambda_1-\lambda_0)\approx 4.6\cdot 10^5$Pa,\newline $\beta=0.28(\mu_1-\mu_0) \approx 4.7\cdot 10^3$Pa,\newline $1000$ pixel, $\#$tetrahedrons $=7644$ (direct problem, offline similar)& $9$min $13$s & $1$h $3$min $52$s & $2$min $12$s\\
\hline
linearized method, \newline $\alpha=0.28(\lambda_1-\lambda_0)\approx 4.6\cdot 10^5$Pa, \newline $\beta=0.28(\mu_1-\mu_0)\approx 4.7\cdot 10^3$ Pa,\newline $2744$ pixel, $\#$tetrahedrons $=37252$ (direct problem, offline similar) & $1$h $4$min & $6$h $40$min & $7$min $39$s\\
\hline
\end{tabular}
\caption{Comparison of the standard (non-linearized) and linearized monotonicity method.}\label{appendix}
\end{center}
\end{table}

\section*{References}
\bibliographystyle{abbrv}
\bibliography{biblio}

\end{document}